\documentclass[final]{siamltex}

\usepackage{amsmath}
\usepackage{amssymb}
\usepackage{amsfonts}
\usepackage{graphicx}
\usepackage{psfrag,bm}
\usepackage{url}
\usepackage{verbatim}
\usepackage{cite}
\usepackage{tikz}
\usepackage{overpic} 
\usetikzlibrary{shapes,arrows,positioning}
\usepackage{psfrag}
\usepackage{subfigure}

\newenvironment{mat}{\left[\begin{array}{ccccccccccccccc}}{\end{array}\right]}
\newcommand\bcm{\begin{mat}}
\newcommand\ecm{\end{mat}}
\newcommand\bigo{\mathcal O}
\newcommand{\Secref}[1]{Section~\ref{#1}} 
\newcommand{\figref}[1]{Figure~\ref{#1}}

\DeclareMathOperator{\cond}{cond}
\DeclareMathOperator{\sign}{sign}
\DeclareMathOperator{\real}{Re}

\usepackage{todonotes}

\newtheorem{problem}{Problem}[section]

\title{Fast computation of Gauss quadrature nodes and weights on the whole real line}
\author{Alex Townsend\thanks{Department of Mathematics, Massachusetts Institute of Technology, 77 Massachusetts Avenue
Cambridge, MA 02139-4307. (ajt@mit.edu, http://math.mit.edu/$\sim$ajt/)}
\and Thomas Trogdon\thanks{Courant Institute of Mathematical Sciences, 251 Mercer Street, New York, 10012-1185. 
(trogdon@cims.nyu.edu, http://www.cims.nyu.edu/$\sim$trogdon/)} \and 
Sheehan Olver\thanks{The University of Sydney, NSW 2006, Australia (Sheehan.Olver@sydney.edu.au, http://www.maths.usyd.edu.au/u/olver/)  }
}

\begin{document}
\maketitle


\begin{abstract} 
A fast and accurate algorithm for the computation of Gauss--Hermite and generalized Gauss--Hermite quadrature nodes and weights is presented. The algorithm is based on 
Newton's method with carefully selected initial guesses for the nodes and a fast evaluation scheme for the associated orthogonal polynomial.  
In the Gauss--Hermite case the initial guesses and evaluation scheme rely on explicit asymptotic formulas. For generalized Gauss--Hermite, the 
initial guesses are furnished by sampling a certain equilibrium measure and the associated polynomial
evaluated via a Riemann--Hilbert reformulation. In both cases the $n$-point quadrature rule
is computed in $\mathcal{O}(n)$ operations to an accuracy that is close to machine precision.  For sufficiently large $n$, some of the quadrature weights have a value less than the smallest positive normalized floating-point number in double precision and we exploit this fact to 
achieve a complexity as low as $\mathcal{O}(\sqrt{n})$.
\end{abstract}

\begin{keywords}
Gauss quadrature, numerical integration, Hermite polynomials, equilibrum measures, Riemann--Hilbert problems
\end{keywords}

\begin{AMS}
65D32, 33C45, 35Q15
\end{AMS}


\section{Introduction}
%
%
%
Numerical quadrature is the approximation of a definite integral of a continuous function $f$ by a weighted linear 
combination of function evaluations, \emph{i.e.}, 
\[
\int_{a}^b f(x) dx \approx \sum_{k = 1}^n w_k f(x_k),  \qquad -\infty \leq a< b \leq \infty,
\]
where $\left\{x_k\right\}_{k=1}^n$ are the {\em nodes} and $\left\{w_k\right\}_{k=1}^n$ are the {\em weights} (indexed so that $x_k<x_{k+1}$). 
An $n$-point quadrature rule of this form is ``Gaussian'' if for some nonnegative weight function, denoted by $w(x)$, the approximation 
\[
\int_{a}^b w(x)f(x) dx \approx \sum_{k = 1}^n w_k f(x_k)
\]
is exact whenever $f$ is a polynomial of degree $\leq 2n-1$. For a fixed weight function this requirement uniquely 
defines a set of quadrature nodes and weights and the resulting integration scheme is called a {\em Gauss quadrature rule}~\cite[Sec.~4.6]{NumericalRecipes}. 

In this paper, we are interested in weight functions of the form $w(x) = e^{-V(x)}$ 
and integrating functions over the whole real line ($a=-\infty$, $b=\infty$). If $V(x)=x^2$, then $w(x)$ is the classic 
Hermite weight. More generally, we are concerned with the so-called Freud weights, where 
$V(x)$ is a real polynomial that grows at infinity.  Of particular interest is the case $V(x) = x^{2m}$, $m\geq 1$. For an integer $n$ and a weight function $w(x)$, there are unique 
sets of Gauss quadrature nodes $\left\{x_k\right\}_{k=1}^n$ and weights $\left\{w_k\right\}_{k=1}^n$, and it is our goal 
to compute these two sets to an accuracy of double precision in $\mathcal{O}(n)$ operations. 

The classic approach for computing Gauss quadrature nodes and weights is the Golub--Welsch algorithm~\cite{Golub_69_01}, 
which requires $\mathcal{O}(n^2)$ operations when one is careful and $\mathcal{O}(n^3)$ operations when one is not\footnote{In many programming languages, for example {\sc Matlab}, 
the structure of symmetric tridiagonal eigenproblems is not automatically detected or exploited.}. 
However, in recent years several fast algorithms have been developed that require only $\mathcal{O}(n)$ operations. Currently, 
the state-of-the-art for classic weight functions is Bogaert's algorithm~\cite{Bogaert_14_01} for Gauss--Legendre ($w(x)=1$, $[a,b] = [-1,1]$), the Hale--Townsend algorithm 
for Gauss--Jacobi~\cite{Hale_13_01} ($w(x)=(1-x)^\alpha(1+x)^{\beta}$, $[a,b]=[-1,1]$), and the Glaser--Lui--Rokhlin algorithm~\cite{Glaser_07_01} for 
Gauss--Laguerre ($w(x)=e^{-x}$, $[a,b]=[0,\infty]$) and Gauss--Hermite ($w(x) = e^{-x^2}$, $[a,b]=[-\infty,\infty]$). In this paper 
we extend the approach in~\cite{Hale_13_01} to a competitive algorithm for computing Gauss--Hermite quadrature nodes and weights. 
Then, we use Riemann--Hilbert (RH) problems to derive an $\mathcal{O}(n)$ algorithm for generalized Gauss--Hermite quadrature
rules. This demostrates, for the first time, that the procedure in~\cite{Hale_13_01} can 
be generalized to nonstandard Gauss quadrature rules.

RH problems are boundary value problems in the complex plane~\cite{AblowitzFokas,Deift_Book,TrogdonThesis}.  
The use of RH problems is necessary for our general approach. The main idea from~\cite{Hale_13_01} that we 
generalize is that asymptotic formulas for orthogonal polynomials and their derivative can be combined with 
initial guesses for the Gauss nodes to derive an effective scheme for computing Gauss quadrature 
nodes and weights.  For classical weights, such expansions are known explicitly, but for other
weight functions we note that accurate approximations can be calculated numerically 
via the solution of a parameter-dependent RH problem~\cite{Deift_Book,TrogdonSORMT} using 
nonlinear steepest descent~\cite{DeiftZhouAMS,TrogdonSONNSD}.

The paper is structured as follows. In the next section we present an overview of our scheme, which 
is central to both the computation of 
Gauss--Hermite nodes and weights and its generalizations. In \Secref{sec:GaussHermite} we describe
how to compute the Gauss--Hermite nodes and weights in $\mathcal{O}(n)$ operations using Newton's 
method together with explicit asymptotic formulas. In \Secref{sec:generalized} we show how these 
ideas can be generalized to weights of the form $e^{-V(x)}$ using equilibrium measures and RH problems. 
Finally, \Secref{sec:interpolation} describes an 
application of these methods to barycentric Lagrange interpolation with an appendix analyzing the weighted 
stability of this interpolation.

\section{Overview of approach}
For a given weight function $w(x)$, our algorithm for computing the corresponding Gauss 
quadrature nodes and weights relies on the standard fact that the nodes are 
precisely the roots of the associated orthogonal polynomial of degree $n$~\cite{Gautschi}. That is, 
if $\phi_0(x),\ldots,\phi_n(x),\ldots,$ is the sequence of orthogonal polynomials that are 
orthogonal with respect to the inner-product 
\[
\langle f,g \rangle = \int_{-\infty}^\infty w(x)f(x)g(x) dx,  
\]
then the Gauss nodes $x_1,\ldots,x_n$ satisfy $\phi_n(x_k) = 0$ for $1\leq k\leq n$. 
This is a powerful observation that transforms the abstract notion of a Gauss quadrature
rule to $n$ tangible rootfinding problems.

We solve each rootfinding problem $\phi_n(x_k) = 0$ with Newton's method, which 
needs three pieces of information: (1) A sufficiently close initial guess for $x_k$, (2) 
An evaluation scheme for $\phi_n$, and (3) An evaluation scheme for $\phi_n'$. As there are 
$n$ rootfinding problems, \emph{i.e.}, $\phi_n(x_k)=0$ for $1\leq k\leq n$, we must solve each one
in just $\mathcal{O}(1)$ operations to achieve an overall complexity of $\mathcal{O}(n)$.
Here is how we achieve (1), (2), and (3): 

\subsection*{(1) Initial guesses}
If $w(x)=e^{-x^2}$, then there are explicit asymptotic expansions that approximate 
the Gauss--Hermite nodes (see~Lemmas~\ref{lem:tricomi} and~\ref{lem:gatteschi}). For $n\geq 200$ these provide sufficiently good  
initial guesses for Newton's method. Unfortunately, for general weights of the form 
$w(x) = e^{-V(x)}$ explicit asymptotic expansions for the Gauss nodes are not available. 
Instead, we use equilibrium measures to furnish initial guesses.  Roughly speaking, the 
equilibrium measure describes the asymptotic density of the \emph{Fekete points}~\cite{SaffPotential}, 
which are the global minimizers of the energy functional
\begin{align*}
E(x_1,\ldots,x_n) = \frac{2}{n(n-1)} \sum_{1 \leq i \neq j \leq n} \log |x_i- x_j|^{-1} + \frac{1}{n} \sum_{i=1}^n V(x_i)
\end{align*}  
as $n \rightarrow \infty$. For $V(x) = x^2$ the Fekete points are the zeros of Hermite polynomials and for general $V(x)$, 
an asymptotic expansion of the zeros of orthogonal polynomials can be derived in terms of the equilibrium measure \cite{DKMVZ}. 

\subsection*{(2) Evaluation of orthogonal polynomial}
If $w(x)=e^{-x^2}$, then the degree $n$ Hermite polynomial can be expressed in terms of a parabolic 
cylinder function, which has a powerful uniform asymptotic expansion involving Airy functions (see~\Secref{sec:fastHermiteEvaluation}). For $n\geq 200$ the resulting asymptotic expansion
is accurate to 14-15 digits and requires only $\mathcal{O}(1)$ operations per evaluation. (Note that we are not able 
to evaluate a Hermite polynomial using the $3$-term recurrence~\cite[(18.9.1)]{NISTHandbook}, since that requires $\mathcal{O}(n)$ operations per 
evaluation.) High-order explicit asymptotic expansions are not available
for orthogonal polynomials associated to generalized Hermite 
weights and instead we use a numerical RH approach (see Section~\ref{sec:RH}).  

\subsection*{(3) Evaluation of the derivative}
In the Gauss--Hermite case, the evaluation of $\phi_n'$ is achieved 
by explicit asymptotic expansions. For more general weights we 
again use a numerical RH approach.

\medskip

Finally, we also need to compute the Gauss quadrature weights. In 
both the Gauss--Hermite and generalized quadrature rules, $w_k$ can be expressed in terms 
of $\phi(x_k)$ and $\phi'(x_k)$; see~\eqref{eq:weightFormula} and~\eqref{weight-formula}.

For sufficiently large $n$, some of the quadrature weights take a value less than the smallest positive normalized 
floating-point number in double precision.  In such circumstances, these quadrature weights and 
corresponding nodes do not contribute to the final quadrature estimate when working in double precision 
(regardless of the function to be integrated). We provide a subsampling scheme so that only the weights and corresponding 
nodes that contribute to the quadrature estimate are computed.  This makes the algorithm far more 
efficient without reducing the accuracy of the resulting quadrature rule.

Though it is not the focus of this paper, when $n$ is small ($n<200$) we recommend using 
Newton's method, as described above, together with polynomial evaluation via the $3$-term 
recurrence~\cite[(18.9.1)]{NISTHandbook}. We observe this approach to be extremely accurate and is 
easily applicable to variable precision computations. For very small $n$, the initial guesses
provided by asymptotics or equilibrium measures will not be sufficient to guarantee convergence of Newton's method. 
At this point, it is reasonable to use the Golub--Welsch algorithm to furnish initial guesses
and Newton's method to improve the accuracy of the final nodes and weights.

\section{Computing Gauss--Hermite quadrature nodes and weights}\label{sec:GaussHermite}
The classical Gauss--Hermite quadrature nodes and weights correspond to the weight function 
$w(x) = e^{-x^2}$, which can be used to approximate the following definite integral: 
\[
 \int_{-\infty}^\infty e^{-x^2} f(x) dx \approx \sum_{k=1}^n w_k f(x_k),
\]
where $f$ is a ``smooth'' function and $\{x_k\}$ and $\{w_k\}$ are the Gauss--Hermite nodes and weights, respectively. The associated sequence of orthogonal polynomials are 
the Hermite polynomials, denoted by $H_0, H_1, \ldots,$ which can be defined via the following $3$-term recurrence 
relation~\cite[(18.9.1)]{NISTHandbook}: 
\[
 H_{n+1}(x) = 2xH_n(x) - 2nH_{n-1}(x), \qquad n\geq 1,\quad x\in\mathbb{R},
\]
where $H_0(x) = 1$ and $H_1(x) = 2x$. The Gauss nodes are the roots of the degree $n$ Hermite polynomial. That is, $H_n(x_k) = 0$ for $1\leq k\leq n$.  

Now, since Hermite polynomials have a reflective symmetry~\cite[(18.6.1)]{NISTHandbook}, \emph{i.e.}, 
$H_n(-x) = (-1)^{n}H_n(x)$, the nodes are symmetrically located on the real line. That is, 
$x_{k} = -x_{n-k+1}$ for $1\leq k \leq n$ and if $n$ is odd, $x_{\lceil n/2\rceil}=0$. 
Therefore, we only need to compute the strictly positive nodes since the others can be obtained by symmetry. 
This will save a factor of roughly $2$ in the computational cost of the final algorithm.

In addition, it is known that the Gauss--Hermite nodes satisfy~\cite[(18.16.16)]{NISTHandbook} 
\[
 -\sqrt{2n+1}<x_1<\cdots<x_n<\sqrt{2n+1},
\]
which means we will only require a fast evaluation scheme for $H_n(x)$ when $0< x< \sqrt{2n+1}$.

\subsection{Initial guesses for Gauss--Hermite quadrature nodes}\label{sec:initialGuesses}
We use two different asymptotic formulas for the strictly positive Gauss--Hermite nodes. For the majority of the nodes 
we use the asymptotic approximations derived by Tricomi~\cite{Tricomi_49_01} and given in the following lemma:
\begin{lemma}[Tricomi~\cite{Tricomi_49_01}]
Let $\tau_{k}$ be the root of the equation 
\begin{equation}
x - \sin x = \frac{(4\lfloor n/2\rfloor-4k+3)\pi}{4\lfloor n/2\rfloor +2\alpha + 2}, \qquad \alpha = {\rm mod}(n,2) - \frac{1}{2}.
\label{eq:Tk}
\end{equation} 
Then, for $n\rightarrow\infty$ and $k\geq 0$ fixed, we have 
\[
x_{k+\lceil n/2\rceil}^2 = \nu \sigma_{k} - \frac{1}{3\nu}\left[\frac{5}{4(1-\sigma_{k})^2} - \frac{1}{1-\sigma_{k}} - \frac{1}{4} \right] + \mathcal{O}(n^{-3}), \quad n\rightarrow\infty,
\]
where $\sigma_{k} = \cos^2(\tau_{k}/2)$ and $\nu = 4\lfloor n/2\rfloor + 2\alpha + 2$. 
\label{lem:tricomi}
\end{lemma} 

In order to compute Tricomi's initial guesses we must first calculate $\tau_k$. To achieve this we 
solve the equation 
in~\eqref{eq:Tk} by several steps of Newton's method with an initial guess of $\pi/2$. If we let $f(x) = x-\sin x$, then 
we note that $f:(0,\pi]\rightarrow(0,\pi]$ is twice differentiable, $f'$ is strictly positive, 
and $f(x)f''(x)>0$ so that convergence of Newton's method to $\tau_k$ is guaranteed regardless of the initial guess.

Tricomi's initial guesses for the nodes are accurate except for a handful near 
$\sqrt{2n+1}$, and for these nodes we use the asymptotic approximations derived by Gatteschi~\cite{Gatteschi_02_01}: 
\begin{lemma}[Gatteschi~\cite{Gatteschi_02_01}]
Let $a_m$ be the $m$th zero of the Airy function ${\rm Ai}(x)$, indexed so that $a_{m+1}<a_m<0$. Then, for $n\rightarrow\infty$ and fixed $k\geq 1$, we have
\[
\begin{aligned}
x_{n-k+1}^2 = \nu + 2^{2/3}a_k\nu^{1/3} &+ \frac{1}{5}2^{4/3} a_k^2 \nu^{-1/3} + \left(\frac{9}{140} - \frac{12}{175}a_k^3\right)\nu^{-1}\\
&+ \left(\frac{16}{1575}a_k + \frac{92}{7875}a_k^4\right)2^{2/3}\nu^{-5/3} \\
&- \left(\frac{15152}{3031875}a_k^5 + \frac{1088}{121275}a_k^2\right)2^{1/3}\nu^{-7/3} + \mathcal{O}(n^{-3}),
\end{aligned}
\]
where $\nu = 4\lfloor n/2\rfloor + 2\alpha + 2$ and $\alpha = {\rm mod}(n,2) - 1/2$. 
\label{lem:gatteschi}
\end{lemma}

In order to compute Gatteschi's initial guesses, the  
zeros of the Airy function are required. We tabulate the first ten roots of 
${\rm Ai}(x)$, while the others are computed with the asymptotic formula~\cite[(9.9.18)]{NISTHandbook} that is observed to be accurate for $m\geq 11$:
\[
a_m \approx -s_m^{2/3}\left(1+\frac{5}{48}s_m^{-2}-\frac{5}{36}s_m^{-4}+\frac{77125}{82944}s_m^{-6} -\frac{108056875}{6967296}s_m^{-8}+\frac{162375596875}{334430208}s_m^{-10}\right),
\]
where $s_m = 3\pi(4m-1)/8$.

In practice, we use Tricomi's initial guesses for $k = 0,\ldots,\lfloor \rho n\rfloor$, where $\rho = 0.4985$, and 
Gatteschi's otherwise. Based on numerical experiments, we have selected $\rho= 0.4985$ because when $n$ is 
large Tricomi's and Gatteschi's initial guesses have roughly the same error for $x_{\lfloor\rho n\rfloor}$. 

\begin{figure}
\centering
\begin{minipage}{.49\textwidth} 
\begin{overpic}[width=\textwidth]{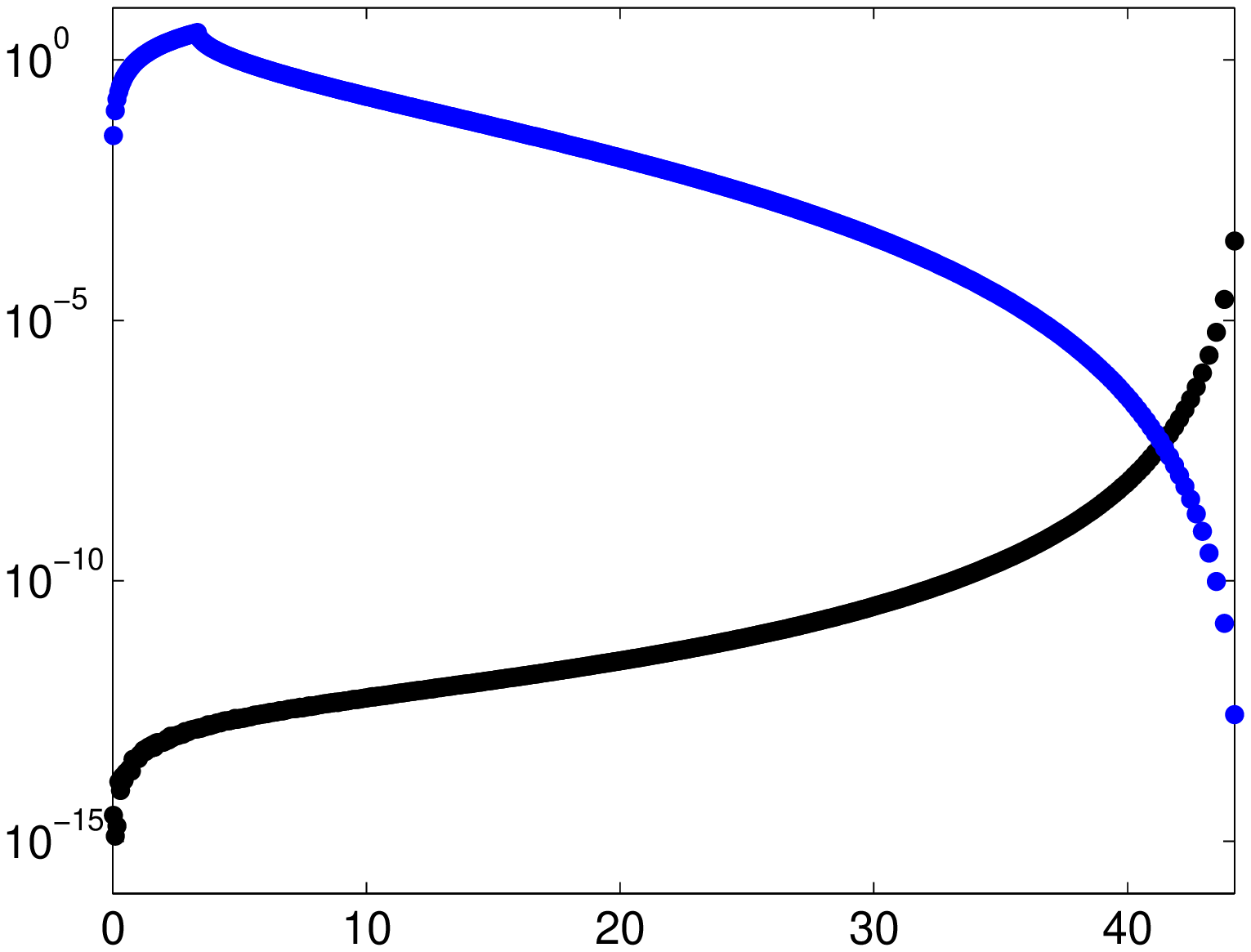}
\put(50,0){$x_k$}
\put(0,20){\rotatebox{90}{Absolute error}}
\put(60,58){\rotatebox{-24}{Lemma~\ref{lem:gatteschi}}}
\put(40,25){\rotatebox{10}{Lemma~\ref{lem:tricomi}}}
\end{overpic}
\end{minipage}
\begin{minipage}{.49\textwidth} 
\begin{overpic}[width=\textwidth]{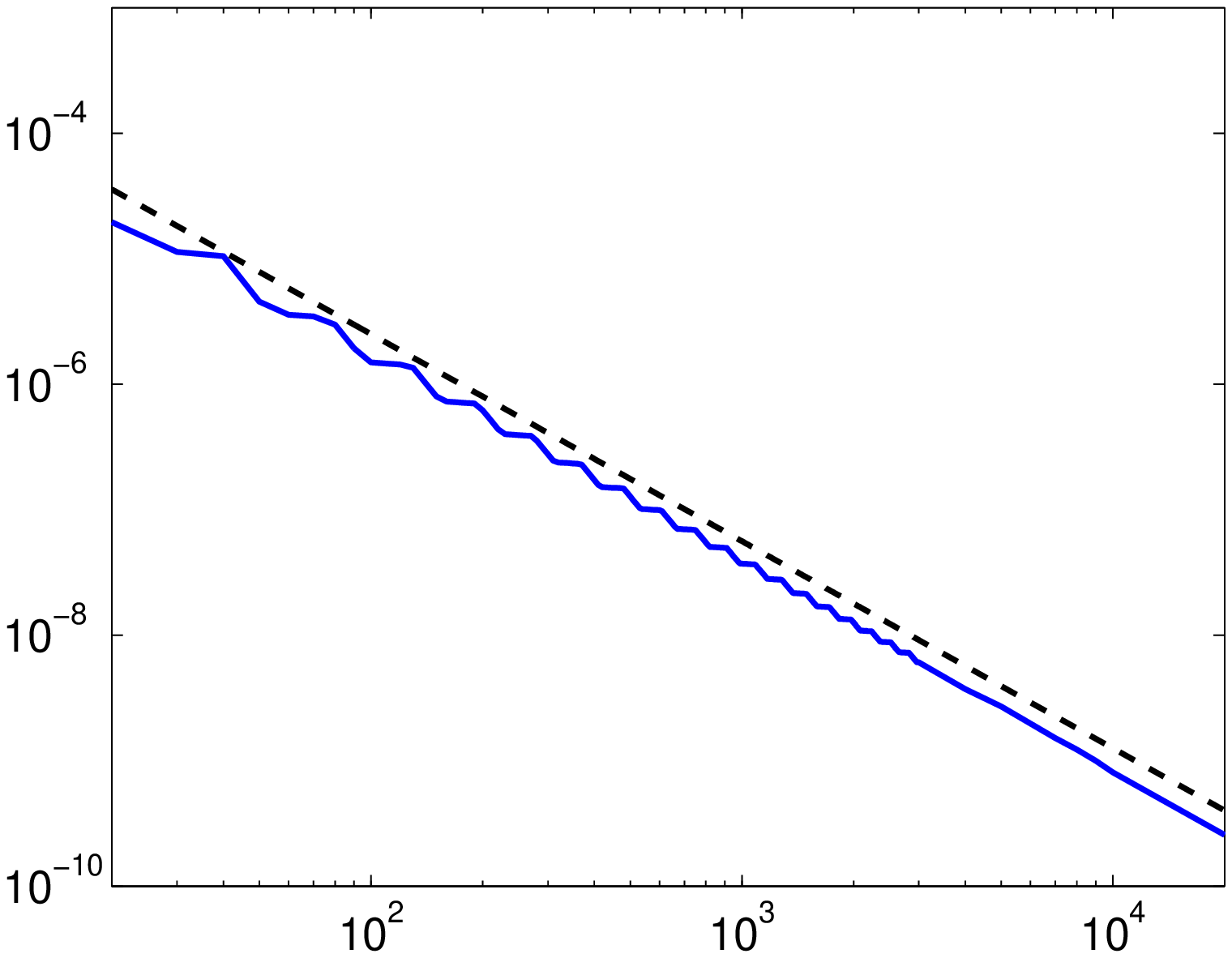}
 \put(50,0){$n$}
\put(0,20){\rotatebox{90}{Absolute error}}
\put(50,40){\rotatebox{-30}{$\mathcal{O}(n^{-1.65})$}}
\end{overpic}
\end{minipage}
\caption{Left: Absolute error of the initial guesses from Lemma~\ref{lem:tricomi} (black) and Lemma~\ref{lem:gatteschi} (blue) for $n=1,\!000$.
Right: The maximum error in the initial guesses for $n\leq 20,\!000$ when using Tricomi's initial guesses for $k=1,\ldots,\lfloor \rho n\rfloor$ and Gatteschi's otherwise, where $\rho = 0.4985$.}
\label{fig:InitialGuessesErrorHermite}
\end{figure}

In~\figref{fig:InitialGuessesErrorHermite} (left) we show the absolute error in Tricomi's and Gatteschi's 
initial guesses for $n = 1,\!000$. It can be seen that Lemma~\ref{lem:tricomi} provides better initial guesses 
except when $x_k\approx \sqrt{2n+1}$. In~\figref{fig:InitialGuessesErrorHermite} (right) we show the absolute error in
the initial guesses for $n \leq 20,\!000$ when Tricomi's initial guesses are used for 
$k = 1,\ldots,\lfloor \rho n\rfloor$ and Gatteschi's otherwise. Interestingly, the observed convergence rate is $\mathcal{O}(n^{-1.65})$, while 
Lemmas~\ref{lem:tricomi} and~\ref{lem:gatteschi} only guarantee a rate of $\mathcal{O}(n^{-1.5})$. 

For $n\geq 6,\!000$, we observe that the initial guesses become so accurate that just one Newton iteration is required 
to compute the Gauss--Hermite nodes to double precision. 

\subsection{Fast evaluation of Hermite polynomials}\label{sec:fastHermiteEvaluation} 
The most powerful asymptotic formulas for Hermite polynomials are based on the asymptotics of the parabolic cylinder function.  
Hermite polynomials satisfy the following relationship~\cite[(18.15.28)]{NISTHandbook}:
\[
H_n(x) = 2^{(\mu^2-1)/4} e^{\mu^2t^2/2} U\left(-\frac{1}{2}\mu^2,\mu t \sqrt{2}\right),
\]
where $U$ is the parabolic cylinder function, $\mu = \sqrt{2n+1}$, and $t = x/\mu$.
Moreover, $U$ has the following asymptotic formula that holds as $\mu\rightarrow\infty$ in the region $-\mu<x\leq\mu$ (equivalently $-1\leq t\leq 1$)~\cite[(12.10.35)]{NISTHandbook}:  
\begin{equation}
U\left(-\frac{1}{2}\mu^2,\mu t\sqrt{2}\right) \sim 2\pi^{\frac{1}{2}} \mu^{\frac{1}{3}} g(\mu) \phi(\zeta)\left({\rm Ai}\left(\mu^{\frac{4}{3}}\zeta\right)\sum_{s=0}^\infty \frac{A_s(\zeta)}{\mu^{4s}} + \frac{{\rm Ai}'\left(\mu^{\frac{4}{3}}\zeta\right)}{\mu^{\frac{8}{3}}}\sum_{s=0}^\infty \frac{B_s(\zeta)}{\mu^{4s}}\right)
\label{eq:cylinderfunction}
\end{equation}
where $\zeta$ satisfies $2/3(-\zeta)^{3/2}=\tfrac{1}{2}\cos^{-1}t - \tfrac{1}{2}t\sqrt{1-t^2}$, and $\phi(\zeta) = (\zeta/(t^2-1))^{1/4}$. Here, 
\[
 g(\mu) = h(\mu) \left(1 + \sum_{s=1}^{\infty} \frac{\phi_s}{(\tfrac{1}{2}\mu^2)^s}\right), \qquad h(\mu) = 2^{-\tfrac{1}{4}\mu^2-\tfrac{1}{4}}e^{-\tfrac{1}{4}\mu^2}\mu^{\tfrac{1}{2}\mu^2-\tfrac{1}{2}}, 
\]
where the coefficients $\phi_s$ are defined by 
\[
 \Gamma(\tfrac{1}{2} + z) \sim \sqrt{2\pi}e^{-z} z^z \sum_{s=0}^\infty \frac{\phi_s}{z^s}.
\]
Moreover, in~\eqref{eq:cylinderfunction} we have, for $t = \cos\theta$,
\[
 A_0(\zeta) = 1, \quad B_0(\zeta) = -(\zeta^6(\cos^3 \theta-6\cos \theta)/24 + 15/144),
\]
and $A_1$, $B_1$, and higher order terms can be calculated from the recurrence~\cite[(12.10.42)]{NISTHandbook}.

\begin{figure}
\centering
\begin{minipage}{.49\textwidth} 
\begin{overpic}[width=\textwidth]{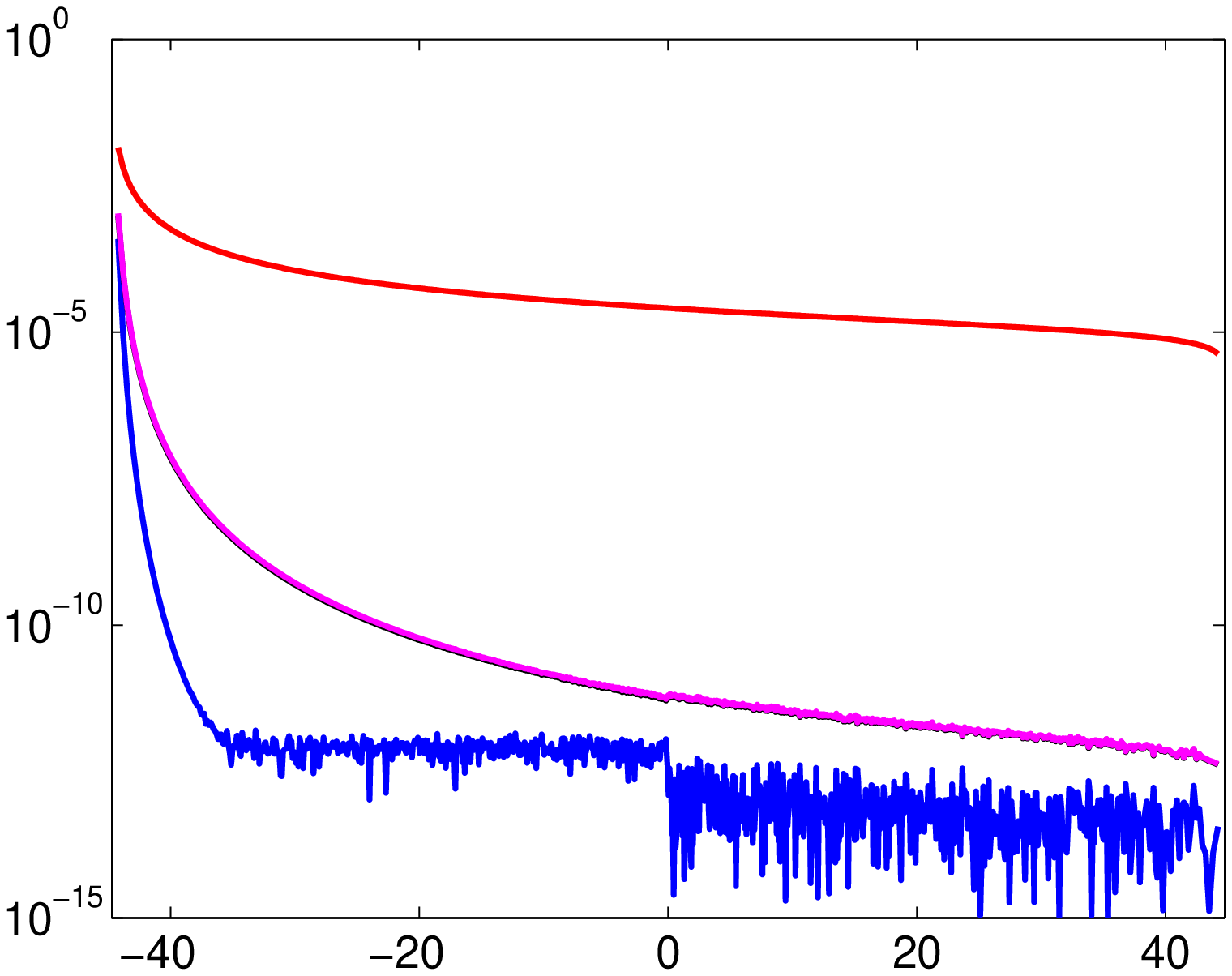}
\put(50,0){$x$}
\put(0,25){\rotatebox{90}{Absolute error}}
\put(32,53){$1$ term}
\put(32,30){\rotatebox{-10}{$2$ \& $3$ terms}}
\put(32,13){$4$ terms}
\end{overpic}
\end{minipage}
\begin{minipage}{.49\textwidth} 
\begin{overpic}[width=\textwidth]{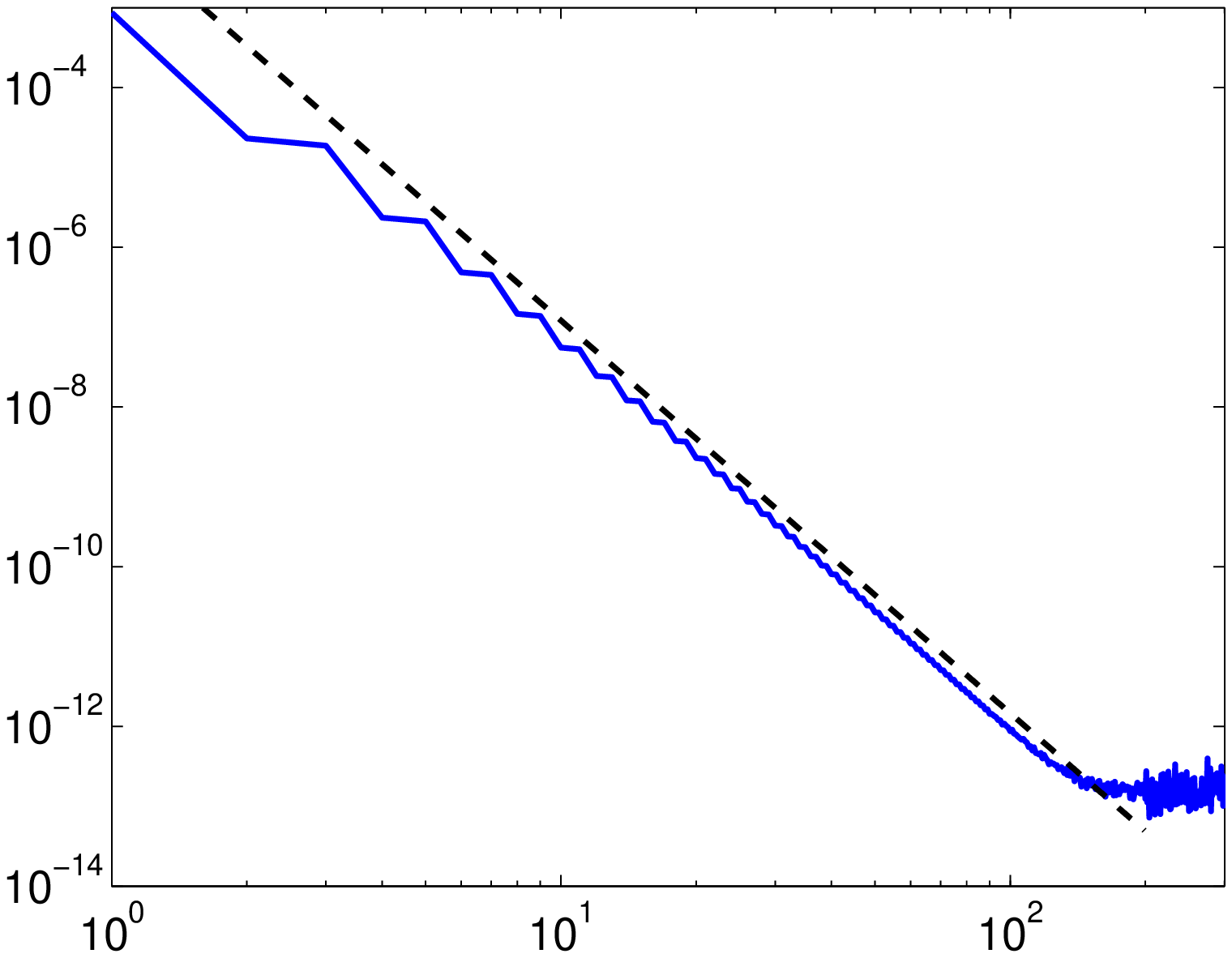}
\put(50,0){$n$}
\put(0,25){\rotatebox{90}{Absolute error}}
\put(55,40){$\mathcal{O}(n^{-59/12})$}
\end{overpic}
\end{minipage}
\caption{Left: Error in the asymptotic formula~\eqref{eq:cylinderfunction} for $\tilde{U}$ with 1 (red), 2 (magenta), 3 (black), and 4 (black) terms when $n=1,\!000$. In our algorithm the expansion is only evaluated at $x>0$ where the asymptotic formula 
is accurate.  Right: Absolute error in the asymptotic formula for $\tilde{U}$ when $1\leq n\leq 250$. 
We observe a convergence rate of $\mathcal{O}(n^{-59/12})$, which is consistent 
with the magnitude of the first neglected term in~\eqref{eq:cylinderfunction}.}
\label{fig:AsymptoticExpansion}
\end{figure}

In practice, we truncate the asymptotic formula in~\eqref{eq:cylinderfunction} after four terms. Of course, 
more (increasingly complicated) terms can be taken, but with just four terms the resulting asymptotic formula is 
accurate for $n\geq 200$ (see~\figref{fig:AsymptoticExpansion}).  Remarkably, despite the rather involved definitions, this asymptotic formula can be evaluated 
to close to $16$~digits of absolute accuracy. 

In order to compute the roots of the Hermite polynomial 
for large $n$, we scale the parabolic cylinder function so that its absolute maximum is 
bounded by $1$. That is, we actually find the roots of 
\begin{equation}
 \tilde{U}\left(-\frac{1}{2}\mu^2,\mu t \sqrt{2}\right) = \frac{2^{\tfrac{1}{4}}U\left(-\frac{1}{2}\mu^2,\mu t \sqrt{2}\right)}{\sqrt{\pi}n^{\tfrac{1}{4}}g(\mu)}, \qquad \left|\tilde{U}\left(-\frac{1}{2}\mu^2,\mu t \sqrt{2}\right)\right|\leq 1.
\label{eq:ScaledCylinderfunction}
\end{equation} 
Scaling in this way is essential for avoiding numerical overflow issues for large $n$. We note that $\tilde{U}$ is closely related to the scaled Hermite polynomial 
described in~\cite[(32)]{Glaser_07_01}.

In~\figref{fig:AsymptoticExpansion} (left) we show the absolute error in the asymptotic formula in~\eqref{eq:cylinderfunction} 
for evaluating $\tilde{U}$ with four terms and $n = 1,\!000$. The asymptotic formula is only evaluated in the region $x>0$ since the 
Gauss--Hermite nodes in $x\leq 0$ can be recovered by symmetry (see Section~\ref{sec:GaussHermite}).  In~\figref{fig:AsymptoticExpansion} (right) we show  
the maximum absolute error of~\eqref{eq:cylinderfunction} in the 
region $x>0$ for $1\leq n \leq 250$.  We observe that the absolute error converges to zero like $\mathcal{O}(n^{-59/12})$. This 
very precise decay rate is expected because (with the scaling in~\eqref{eq:ScaledCylinderfunction}) the first neglected term (fifth term) 
of the asymptotic formula is of magnitude $\mathcal{O}( {\rm Ai}(\mu^{\frac{4}{3}}\zeta)n^{\frac{1}{12}}\mu^{-8}) = \mathcal{O}(n^{-59/12})$.
\subsection{Fast evaluation of the derivative of Hermite polynomials}
To evaluate $H_n'(x)$, we have two options. One option is to use the following differentiation relationship~\cite[(18.9.25)]{NISTHandbook}:
\[
 H_n'(x) = 2n H_{n-1}(x), \qquad n\geq 1, 
\]
and then to employ~\eqref{eq:cylinderfunction} corresponding to $H_{n-1}$ rather than $H_n$. Another option, 
and the one we employ, is to use the asymptotic formula for $U'$ given by~\cite[(12.10.36)]{NISTHandbook}
\begin{equation}
U'\left(-\frac{1}{2}\mu^2,\mu t\sqrt{2}\right) \sim \frac{(2\pi)^{\frac{1}{2}} \mu^{\frac{2}{3}} g(\mu)}{\phi(\zeta)}\left(\frac{{\rm Ai}\left(\mu^{\frac{4}{3}}\zeta\right)}{\mu^{\frac{4}{3}}}\sum_{s=0}^\infty \frac{C_s(\zeta)}{\mu^{4s}} + {\rm Ai}'\left(\mu^{\frac{4}{3}}\zeta\right)\sum_{s=0}^\infty \frac{D_s(\zeta)}{\mu^{4s}}\right),
\label{eq:derivativecylinderfunction}
\end{equation}
where we have, for $t = \cos\theta$,
\[
 C_0(\zeta) = \frac{2}{3}\left(\zeta^6(\cos^3\theta + 6\cos\theta)/24 - 7/48\right)\zeta^{-\frac{2}{3}}, \qquad D_0(\zeta) = 1,
\]
and higher order terms can be obtained from the recurrence~\cite[(12.10.44)]{NISTHandbook}.  We prefer the latter approach 
because~\eqref{eq:derivativecylinderfunction} contains exactly the same Airy functions as~\eqref{eq:cylinderfunction} 
and hence, the expensive (but $\mathcal O(1)$) special function evaluations can be reused.  

\subsection{Newton's method for Gauss--Hermite nodes} 
We now have all the ingredients to compute Gauss--Hermite nodes using Newton's method. However, we do not perform 
Newton's method in the usual $x$-variable, but instead the $\theta$-variable, where 
\[
 t = x/\mu, \qquad t = \cos \theta. 
\]
This improves the accuracy of the final nodes, particularly those close to $x\approx \mu$, \emph{i.e.}, $t\approx 1$, because it does not require an evaluation of $\cos^{-1}(\,\cdot\,)$ per 
iteration, which is sensitive to small perturbations in arguments close to $1$. Therefore, we take the initial guesses from~\Secref{sec:initialGuesses} 
and perform a change of variables to obtain initial guesses in the $\theta$-variable. We then proceed with Newton's method performed in this variable. 
For the $k$th Gauss--Hermite node one Newton step takes the form: 
\[
 \theta_k^{\rm new} = \theta_k^{\rm old} + \frac{\tilde{U}\left(-\frac{1}{2}\mu^2,\sqrt{2}\mu \cos\theta_k^{\rm old} \right)}{\sqrt{2}\mu \tilde{U}'\left(-\frac{1}{2}\mu^2,\sqrt{2}\mu \cos\theta_k^{\rm old} \right)\sin\theta_k^{\rm old}}. 
\]
If the update, $|\theta_k^{\rm new} - \theta_k^{\rm old}|$, is sufficiently small then Newton's method is terminated and the corresponding 
Gauss--Hermite node is calculated via $x_k = \mu\cos\theta_k^{\rm new}$.    In practice, we use the 
same number of iterations for every Gauss node so that Newton's method can 
be vectorized for a slightly improved computational efficiency. 

\subsection{Computing Gauss--Hermite weights} 
Once the Gauss--Hermite nodes have been computed, the Gauss--Hermite weights immediately follow by the simple formula~\cite[(40)]{Glaser_07_01}: 
\[
w_k = 2e^{-x_k^2} / \tilde{H}'^2_n(x_k),
\]
where $\tilde{H}_n$ is the Hermite polynomial scaled so that $|\tilde{H}_n(x)|\leq1$ for $x>0$. In terms of parabolic cylinder functions
this results in the following formula: 
\begin{equation}
w_k = \frac{Ce^{-x_k^2}}{\left(\tilde{U}'\left(-\frac{1}{2}\mu^2,\mu t_k \sqrt{2}\right)\right)^2}, \qquad t_k = x_k/\mu,
\label{eq:weightFormula}
\end{equation}
where $C$ is a constant so that $\sum_k w_k = \sqrt{\pi}$.

\subsection{Subsampling}\label{subsec:GHsubsampling} 
The Gauss--Hermite quadrature weights that correspond to nodes far away from $0$ are usually 
very small in magnitude. So much so, that a significant proportion of the quadrature weights 
are less than {\tt realmin}, \emph{i.e.}, $2^{-1022} \approx 2.23\times 10^{-308}$, which is the smallest normalized positive 
floating-point number in double precision. Thus, for any quadrature rule employed in double precision these weights 
will never contribute to the final approximation of the integral.  

We have observed that the only quadrature weights that are larger than {\tt realmin} are $w_{\lfloor n/2\rfloor-M+1},\ldots,w_{\lfloor n/2\rfloor+M}$ when $n$ is even and  $w_{\lfloor n/2\rfloor-M},\ldots,w_{\lfloor n/2\rfloor+M}$ when $n$ is odd, where $M=\lceil 12.5 n^{1/2}\rceil$.  
Therefore, we only run Newton's method with initial guesses for this subset of the weights. Since 
$M=\mathcal{O}(n^{1/2})$ the resulting algorithm has a reduced complexity 
of $\mathcal{O}(n^{1/2})$ operations without sacrificing the accuracy of the resulting quadrature rule. In Table~\ref{tab:subsampling}
the execution time in seconds is given for computing the Gauss--Hermite nodes and weights with and without subsampling. 
It is observed that a significant saving can be achieved by not computing weights that have a magnitude less than {\tt realmin}.  
Similar computational savings can be easily achieved by the REC and GLR algorithms (see Section~\ref{sec:numericalResults}). 
We believe it is much harder to avoid the computation of a subset of the nodes and weights in the Golub--Welsch algorithm. 

\begin{table}
\centering
\begin{tabular}{cccc}
$n$ & $M$ & No subsampling &  Subsampling \cr 
\hline 
$10$ & $5$ & $0.004951$ & $0.005011$ \cr
$100$ & $50$ & $0.005663$ & $0.005432$ \cr
$1,\!000$ & $361$ & $0.011345$ & $0.009602$\cr
$10,\!000$ & $1,\!212$ & $0.053017$ & $0.027845$\cr
$100,\!000$ & $3,\!848$ & $0.380275$ & $0.049851$\cr
$1,\!000,\!000$ & $12,\!156$ & $3.635814$ & $0.121961$\cr
\hline
& $\mathcal{O}(n^{1/2})$ & $\mathcal{O}(n)$ & $\mathcal{O}(n^{1/2})$\cr
\end{tabular}
\caption{Execution time in seconds for computing Gauss--Hermite nodes and weights with and without subsampling. For the $n$-point Gauss--Hermite quadrature rule only $\mathcal{O}(\sqrt{n})$ weights contribute to the final quadrature approximation in double precision.  Hence, a significant proportion of the computation of Gauss rules on the whole real line can be saved.}
\label{tab:subsampling}
\end{table} 

\subsection{Numerical results}\label{sec:numericalResults} 
In this section we compare the algorithm described in this section 
based on asymptotic formulas (ASY for short) against three other methods for computing Gauss--Hermite 
quadrature nodes and weights, which we refer to using 
the acronyms:
\begin{description}
 \item {REC:} This recurrence-based algorithm performs Newton's method with orthogonal polynomial evaluation using a $3$-term recurrence, requiring $\mathcal{O}(n^2)$ operations;  
 \item {GLR:} The Glaser--Lui--Rokhlin algorithm solves the associated 2nd-order ordinary differential equation~\cite[Table 18.8.1]{NISTHandbook} using a predictor-corrector-like marching scheme, requiring $\mathcal{O}(n)$ operations~\cite{Glaser_07_01};
 \item {GW:} The Golub--Welsch algorithm solves for the eigenvalues (nodes) and eigenvectors (related to the weights) of 
 the associated Jacobi matrix. The implementation we use for GW here requires $\mathcal{O}(n^3)$ operations as we do not 
 exploit the symmetric tridiagonal structure of the Jacobi matrix~\cite{Golub_69_01}. 
\end{description}
These algorithms have been implemented in MATLAB and the numerical comparisons are performed in that language. 
For the GLR algorithm we use the MATLAB implementation in the {\tt hermpts} command in Chebfun~\cite{Chebfun}. 
The algorithm for the ASY method has also been implemented in the {\tt gausshermite} command in 
the FastGaussQuadrature package~\cite{FastGaussQuadrature}, which is written in the Julia language~\cite{Julia}. 

\begin{figure}
\centering
\begin{minipage}{.49\textwidth} 
\begin{overpic}[width=\textwidth]{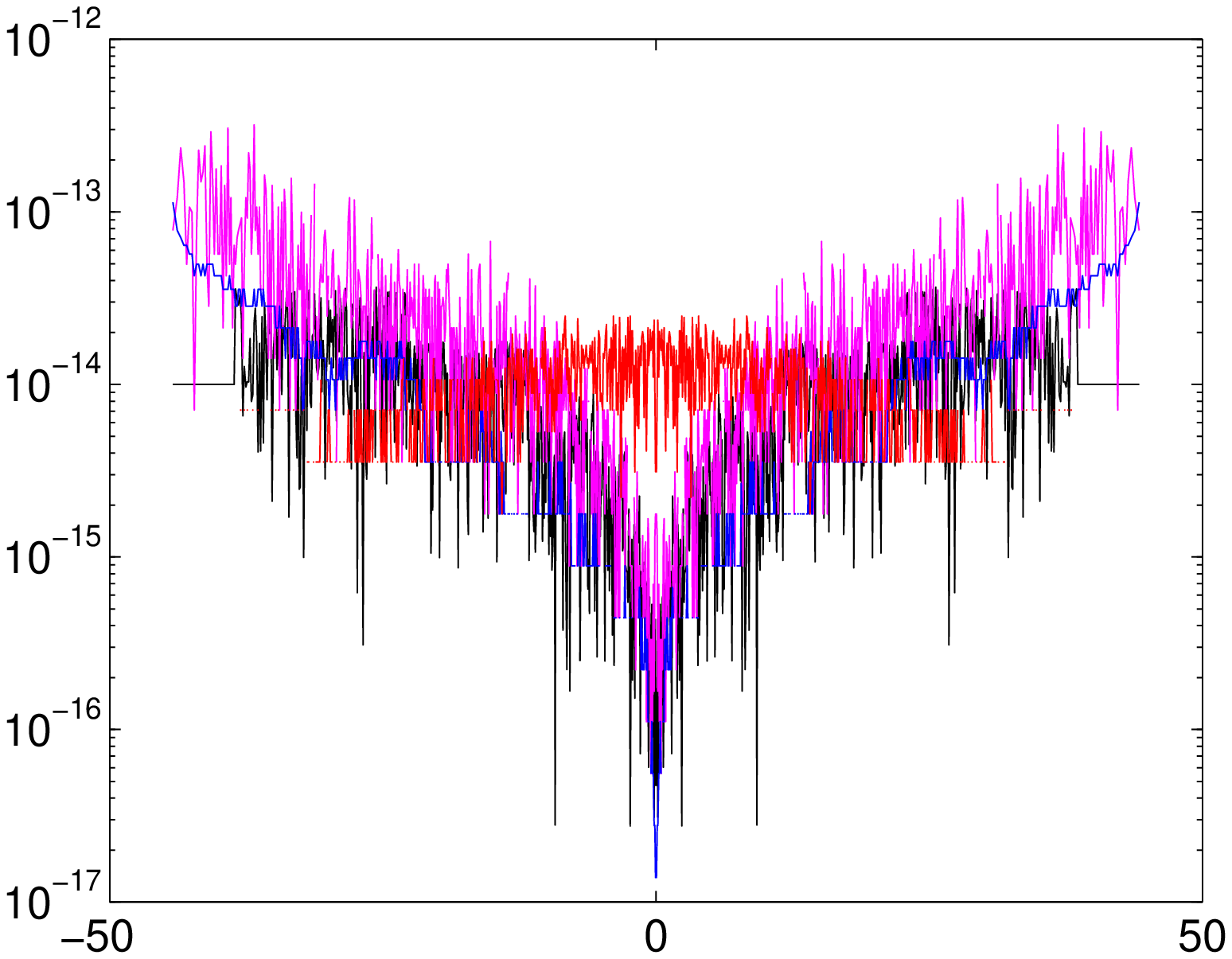}
\put(50,0){$x_k$}
\put(0,25){\rotatebox{90}{Absolute error}}\end{overpic}
\end{minipage}
\begin{minipage}{.49\textwidth} 
\begin{overpic}[width=\textwidth]{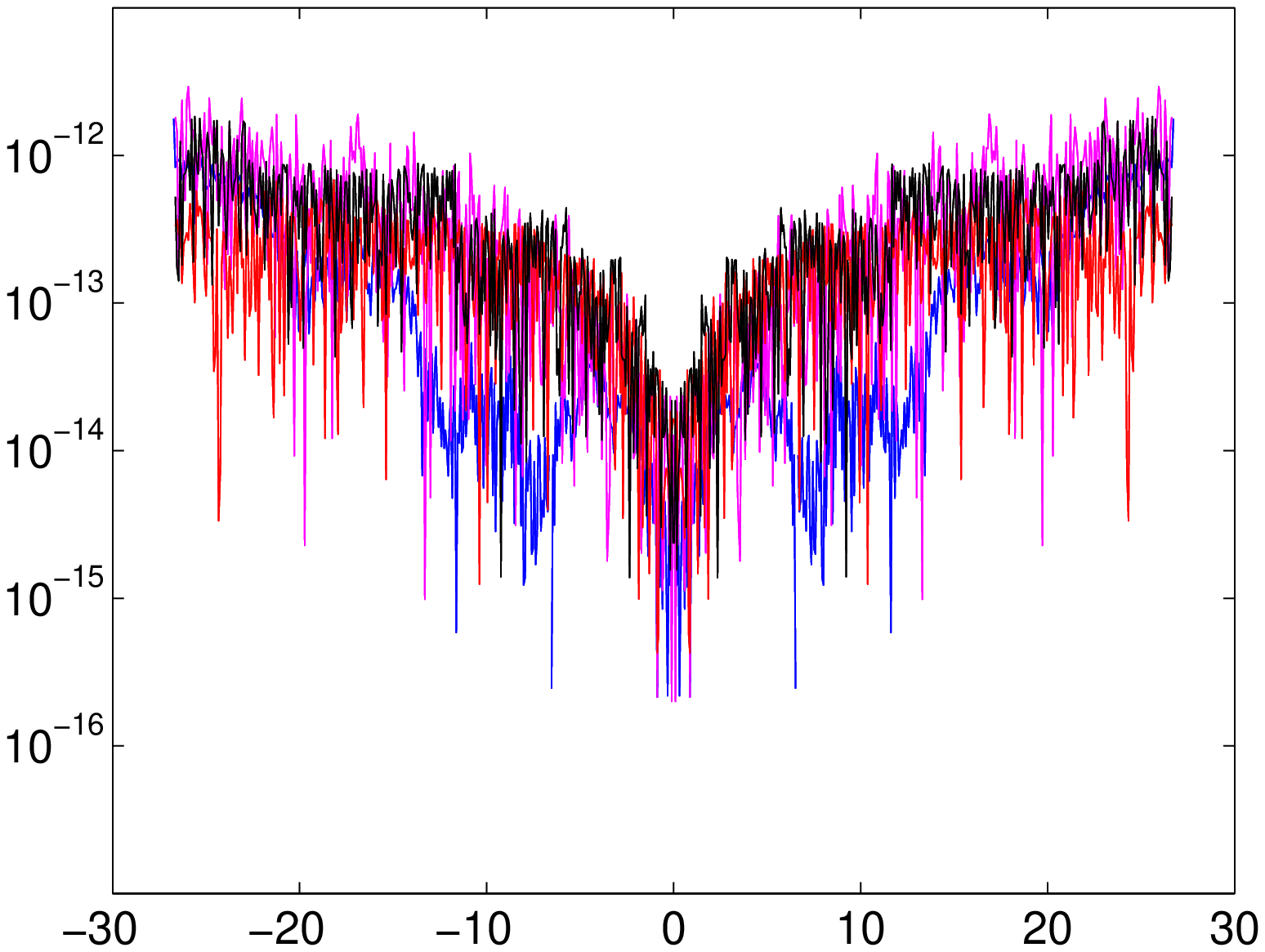}
\put(50,0){$w_k$}
\put(0,25){\rotatebox{90}{Relative error}}
\end{overpic}
\end{minipage}
\caption{Left: Absolute error in the computed Gauss--Hermite nodes, \emph{i.e.}, $|x_k - x_k^{quad}|$, where $x_k^{quad}$ are 
the nodes obtained using quadruple precision, for ASY (red), REC (black), GLR (blue), and GW (magenta).  Right: Relative error in the computed Gauss--Hermite weights, \emph{i.e.}, 
$|w_k - w_k^{quad}|/|w_k^{quad}|$, where $w_k^{quad}$ are the weights obtained using quadruple precision for ASY (red), REC (black), GLR (blue), and GW (magenta). Weights corresponding
to nodes with absolute magnitude larger than $27$ are less than the smallest 
positive normalized floating-point number in double precision and underflow.}
\label{fig:HermiteErrors}
\end{figure}

In Figure~\ref{fig:HermiteErrors} we show the absolute errors in the computed Gauss--Hermite nodes $|x_k - x_k^{quad}|$ (left) and 
the relative error in the weights $|w_k - w_k^{quad}|/|w_k^{quad}|$ (right), where $x_k^{quad}$ and $w_k^{quad}$ are 
the nodes and weights computed using REC with quadruple precision. The ASY computed nodes are less accuracy near $x=0$ but 
more accurate near $x=\sqrt{2n+1}$, which is caused by the 
asymptotic formula~\eqref{eq:cylinderfunction} written in terms of the $\theta$-variable. More accurate Gauss--Hermite nodes could be
obtained by a hybrid between the nodes computed by GLR and ASY. 

\begin{figure}
\centering
\begin{overpic}[width=.49\textwidth]{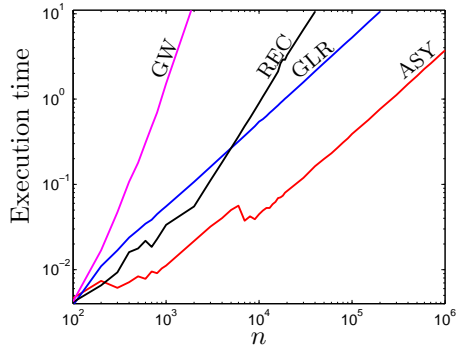}
\put(50,0){$n$}
\put(0,25){\rotatebox{90}{Execution time}}
\put(29,55){\rotatebox{68}{\footnotesize{GW}}}
\put(51,55){\rotatebox{57}{\footnotesize{REC}}}
\put(58,55){\rotatebox{40}{\footnotesize{GLR}}}
\put(80,55){\rotatebox{40}{\footnotesize{ASY}}}
\end{overpic}
\caption{Comparison of computational times for Gauss--Hermite nodes and weights: ASY (red), REC (black), GLR (blue), and GW (magenta).}
\label{fig:HermiteTimings}
\end{figure}

In Figure~\ref{fig:HermiteTimings} we compare the computational timings for the four methods. It can be seen 
that ASY is about $10$ times faster than GLR for large $n$ and ASY takes the lowest execution time out of the 
four methods when $n\geq 200$. The kink in the ASY timings at $n=6,\!000$ is caused by one fewer Newton iteration required 
for convergence when $n\geq 6,\!000$. For $n=1,\!000,\!000$, ASY requires $3.63$ seconds. Subsampling can improve the 
computational speed of ASY, REC, and GLR (see Section~\ref{subsec:GHsubsampling}), but the comparisons between 
these three methods will stay the same. 

\section{Computing generalized Gauss--Hermite quadrature nodes and weights}\label{sec:generalized}

The generalized Gauss--Hermite quadrature nodes and weights correspond to the weight function 
$w(x) = e^{-V(x)}$, where $V(x) = x^{2m} + \bigo(x^{2m-1})$ is a monic polynomial of degree $2m$ with real coefficients. As before, for an integer $n$, we denote the set of generalized Gauss--Hermite nodes and weights as $\{x_k\}_{k=1}^n$ and $\{w_k\}_{k=1}^n$, respectively. 

\subsection{Initial guesses for generalized Gauss--Hermite quadrature nodes}\label{sec:equilibrium}

Given a polynomial $V(x) = x^{2m}+ \mathcal O(x^{2m-1})$ with real coefficients, the equilibrium measure $d\mu_{V,n}(x) = \psi_{V,n}(x) dx$ is the unique minimizer of the functional~\cite{Deift_Book}
\begin{align*}
H(\mu) = \iint \log |x-y|^{-1} d\mu(x)d\mu(y) + \frac{1}{n}\int V(x n^{1/(2m)}) d\mu(x),
\end{align*}
among Borel probability measures $\mu$ on $\mathbb R$, \emph{i.e.}, $\psi_{V,n}(x) dx = \mathrm{argmin}_\mu H(\mu)$.  Without loss of generality, we assume $V(0) = V'(0) = 0$ (otherwise perform a change of variables), 
so that the support of the equilibrium measure $[a,b] \equiv [a_n,b_n]$ converges as $n\rightarrow \infty$ to an interval of the form $[-\alpha,\alpha]$, $\alpha > 0$.
The zeros $\{\tilde x_{k}\}_{k=1}^n$ of the $n$th-order polynomial with varying weight $e^{-V(x n^{1/(2m)})}dx$ are distributed according to the equilibrium measure in the sense that the normalized counting measure converges in distribution~\cite{Deift_Book} (see Figure~\ref{EightHisto}). That is,
\begin{align}\label{dist-conv}
\lim_{n\rightarrow\infty} \left(\frac{1}{n} \int_{a}^t \sum_{k=1}^n \delta_{\tilde x_{k}}(x) dx - \int_a^t d\mu_{V,n}(x) \right) = 0, \quad t \in [a,b],
\end{align}
where $\delta_{\tilde x_{k}}(x)$ is the delta function centered at $\tilde{x}_k$. Therefore, the equilibrium measure can be used to determine the approximate location of the zeros and hence, 
used to furnish initial guesses for Gauss nodes.  It is also known that $\tilde{x}_k \in (a,b)$~\cite[Prop.~3.42]{Deift_Book}.  We consider polynomials with respect to the varying weight $e^{-V(x n^{1/(2m)})}dx$ because it is easily seen that the Gauss quadrature nodes with respect to $e^{-V(x)}dx$ satisfy $ x_{k} = \tilde x_{k} n^{1/(2m)}$ for $1\leq k\leq n$.  Furthermore, the weights $\{w_k\}_{k=1}^n$ for $e^{-V(x)}dx$ satisfy $w_k = \tilde w_k n^{1/(2m)}$, where $\{\tilde w_k\}_{k=1}^n$ are the quadrature weights corresponding to the varying 
weight function $e^{-n V(xn^{1/(2m)})}$.

Define the function $F_n:[a,b]\rightarrow \mathbb{R}$ as 
\[
F_n(x) = \begin{cases} 0, &  x<a,\\
\int_{a}^x \psi_{V,n}(t) dt, & a\leq x \leq b,\\
1, & x > b,\end{cases}
\]
and the inverse function $G_n:(0,1)\rightarrow \mathbb{R}$ as $G_n(y) =  \inf_{x \in \mathbb R}\{F_n(x) \geq y\} $. Then, for $y \in (0,1)$ the value of $G_n(y)$ lies in the support of $\psi_{V,n}(x)$ and $F_n(G_n(y)) = y$.  The following asymptotic formula for the generalized Gauss--Hermite nodes is given in~\cite{DKMVZ}:
\begin{align*}
\left|\tilde x_{k} - G_n\left( \frac{2k-1}{2n} + \frac{1}{2 \pi n}\sin^{-1}\left(G_n(k/n)\right) \right) \right| < \frac{C}{n^2\left(\frac{k}{n}(1-\frac{k}{n})\right)^{4/3}},
\end{align*}
where $C$ is a constant.  Therefore, we can furnish Newton's method with the following initial guesses: \newcommand{\gu}{\mathrm{guess}}
\begin{align*}
\tilde x_k^\gu = G_n\left( \frac{2k-1}{2n} + \frac{1}{2 \pi n}\sin^{-1}\left(G_n(k/n)\right)\right).
\end{align*}

\subsubsection{Computing $\bm{\psi}_{\mathbf{V,n}}$ and $\bm{G_n}$}

While the definition of an equilibrium measure is stated as an optimization problem over measures,  in the case of smooth $V$ the support of 
the equilibrium measure $\psi_{V,n}$ is a single interval $[a,b]$ for sufficiently large $n$. The problem reduces 
to a simpler optimization over the two parameters $a$ and $b$, which can be efficiently solved by Newton's method~\cite{OlverEM}.  
In particular, when $V(x)$ is a polynomial, the equilibrium measure has the form $\psi_{V,n}(x) = \sqrt{(b-x)(x-a)}p(x)$ for some 
polynomial $p(x)$ of degree $2m-2$~\cite[p.~175]{Deift_Book}.  Once $[a,b]$ is found,  $p(x)$ can be calculated in a Chebyshev expansion of the second kind so that
\begin{align*}
\psi_{V,n}(x) = \sqrt{(b-x)(x-a)} \sum_{j=0}^{2m-2} \beta_j U_j\left(M(x)\right), \qquad M(x) = \frac{2x - b-a}{b-a},
\end{align*} 
where $U_j$ is the degree $j$ Chebyshev polynomial of the second kind and the coefficients $\beta_j$ are determined by
	$$V'(x) =  \pi (b-a) \sum_{j=0}^{2m-2} \beta_j T_{j+1}\left(M(x)\right).$$
The interval $[a,b]$ is selected so that the zero-th Chebyshev coefficient in the above expansion vanishes.
The function 
$\psi_{V,n}(x)$ is also important for evaluating the associated orthogonal polynomials using a 
RH reformulation (see Section~\ref{sec:RH}).

Once $\psi_{V,n}(x)$ has been computed, we compute $G_n$ by first calculating $F_n$. For any $x\in[a,b]$, $F_n(x)$ 
is defined by an indefinite integral of $\psi_{V,n}(x)$ so we rewrite the expansion for $\psi_{V,n}(x)$ in terms 
of the Chebyshev basis of the first kind using the recurrence~\cite[(18.9.10)]{NISTHandbook}:
\begin{align*}
\psi_{V,n}\left(M^{-1}(x) \right) &= \frac{b-a}{2} \frac{1}{\sqrt{1-x^2}} \sum_{j=0}^{2m} \frac{\beta_{j-2}-\beta_j}{2}T_j(x),
\end{align*}
where $\beta_{-2}=\beta_{-1}=\beta_{2m-1}=\beta_{2m}=0$ and $T_j(x) = \cos(j\cos^{-1}x)$ is the degree $j$ Chebyshev polynomial of the first kind. (A factor of $(1-x^2)$ appears in the recurrence~\cite[(18.9.10)]{NISTHandbook}, which results in the $(1-x^2)^{-1/2}$ term above.)  This 
allows the indefinite integral to be easily calculated since $T_j$ satisfies the following relation: 
\begin{align}\label{eq:cumsumChebT}
\int_{-1}^y T_j(t) \frac{dt}{\sqrt{1-t^2}} = \begin{cases} -\frac{\sin (j \cos^{-1} y)}{j}, & j >0,\\
\pi - \cos^{-1} y, & j = 0,
\end{cases}\qquad y\in[-1,1].
\end{align}

Therefore, for any $x\in[a,b]$ we have
\[
F_n(x) = \int_{a}^x \psi_{V,n}(t) dt = \frac{(b-a)^2}{4} \int_{-1}^{M(x)} \sum_{j=0}^{2m} \frac{\beta_{j-2}-\beta_j}{2}T_j(t) \frac{dt}{\sqrt{1-t^2}},
\]
where the indefinite integral can be calculated by~\eqref{eq:cumsumChebT}. 
The function $G_n(x)$ can now be computed with Newton's method applied to $F_n(x)$ since evaluation of $F_n(x)$ and its derivative can 
be computed efficiently and accurately.

\subsection{Fast evaulation of generalized Hermite polynomials}\label{sec:RH}

Evaluating orthogonal polynomials associated to the weight function $e^{-V(x)}dx$ with $V(x) = x^{2m} + \mathcal O (x^{2m-1})$ can be achieved using RH techniques. We do not present the full method here as further details can be found in~\cite{SORHFramework, TrogdonSORMT, TrogdonSOGaussQuad}.  
At its very essence the method solves the following RH problem:
\begin{problem}\label{RHP}
Find $\Phi: \mathbb C \setminus \mathbb R \rightarrow \mathbb C^{2\times 2}$ such that $\Phi$ is analytic in $\mathbb C \setminus \mathbb R$ and satisfies
\begin{align*}
\lim_{\epsilon \rightarrow 0^+} \Phi(x + i \epsilon) &= \lim_{\epsilon \rightarrow 0^+} \Phi(x - i \epsilon) \begin{mat} 1 & e^{-V(x n^{1/(2m)})} \\ 0 & 1 \end{mat},\\
\lim_{z \rightarrow \infty} \Phi(z)& \begin{mat} z^{-n} & 0 \\ 0 & z^{n} \end{mat} = \begin{mat} 1 & 0 \\ 0 & 1 \end{mat}.
\end{align*}
\end{problem}

\newcommand{\CC}{\mathcal C}

Remarkably, it can be shown that $\Phi(z)$ takes the following form~\cite{Deift_Book}:
\begin{align*}
\Phi(z) &=\begin{mat}
  \pi_{n}(z)& \CC_{\mathbb R}{[\pi_{n}(\,\cdot\,) e^{-V(\,\cdot\, n^{1/(2m)})}](z)}  \\
   -2 \pi i \gamma_{n-1} \pi_{n-1}(z) &  - 2 \pi i\gamma_{n-1} \CC_{\mathbb R}[\pi_{n-1}(\,\cdot\,) e^{-V(\,\cdot\, n^{1/2m})}](z)
\end{mat},\\
\CC_{\mathbb R}[f](z) &= \frac{1}{2\pi i} \int_{\mathbb R} \frac{f(s)}{s-z} ds,
\end{align*}
where $\pi_j(x)$ is the degree $j$ monic orthogonal polynomial associated to the weight function $e^{-V(x n^{1/(2m)})}$ and $\gamma_n$ is a normalization constant, \emph{i.e.}, 
\begin{align*}
\gamma_n = \left[ \int_{\mathbb R} \pi_{n}^2(x) e^{-V(x n^{1/(2m)})} dx\right]^{-1}.
\end{align*}

Initially, $\Phi(z)$ has growth in the first column and decay in the second at infinity and this is not immediately tractable for the numerical method of \cite{SORHFramework}.  Instead, define
\begin{align*}
g(z) = \int \log(z-x) \psi_{V,n}(x) dx.
\end{align*}
so that  $g$ has a branch cut on and to the left of the support of $\psi_{V,n}$.  It then follows that
\begin{align*}
e^{n g(z)} = z^n + \bigo(z^{n-1}) \quad \text{as} \quad z \rightarrow \infty
\end{align*}
is analytic away from the support of $\psi_{V,n}$.  The function $g(z)$ has many important properties for the asymptotic analysis of Problem~\ref{RHP}.  Two properties that are important for computation are: 
\begin{itemize}
\item $e^{n g(z)}$ captures the growth of $\pi_n(z)$ at infinity,
\item $e^{n g(z)}$ captures the oscillatory behavior of $\pi_n$ on the support of $\psi_{V,n}$.
\end{itemize}

The method in \cite{TrogdonSORMT} provides a numerical solution to the RH problem and 
returns a function $U: \Gamma(n) \rightarrow \mathbb C^{2\times 2}$, where $U$ and $\Gamma(n)$ satisfy
\begin{align*}
U|_{\Gamma_{n,j}}(z) &= \sum_{N=0}^{m(n)} \alpha_{j,N} T_N \left( \frac{z-d_{n,j}}{c_{n,j}}  \right),\\
\Gamma_{n,j} &= c_{n,j} [-1,1] + d_{n,j},\\
\Gamma(n) &= \Gamma_{n,1} \cup \cdots \cup \Gamma_{n,N}.
\end{align*}
Here, the constants $c_{n,j}$ and $d_{n,j}$ are determined by the deformation of the RH problem and $\alpha_{j,k}$ are matrix-valued constants.  There exists a matrix-valued function $N(z):\mathbb{C}\rightarrow\mathbb{C}^{2\times 2}$ 
such that~\cite{TrogdonSORMT}
\begin{align}\label{e:Soln}
\Phi(z) &\approx \tilde \Phi(z) = L_n (\CC_{\Gamma(n)}[U](z)+ I)N(z) L_n^{-1} \begin{mat} e^{n g(z)} & 0 \\ 0 & e^{-ng(z)} \end{mat},\\
\CC_{\Gamma(n)}[U](z) &= \frac{1}{2\pi i} \int_{\Gamma(n)} \frac{U(s)}{s-z} ds,\notag
\end{align}
where $I$ is the $2\times 2$ identity matrix and $L_n$ is a constant diagonal matrix
\begin{align*}
L_n = \begin{mat} e^{-n\ell/2}  & 0 \\ 0 & e^{n\ell/2} \end{mat}.
\end{align*}
The constant $\ell$ is determined from $g$.  Typically, the error made in this approximation is on the order of machine precision.   As described in \cite{TrogdonSORMT}, $m(n)$ is bounded as a function of $n$.  The precise form of $N(z)$ can be deduced from \cite[Section 4]{TrogdonSORMT}.

Evaluating $\mathcal C_{\Gamma(n)} [U](z)$ requires $\mathcal O(1)$ operations because $m(n)$ is bounded as a function of $n$.  The other factors in \eqref{e:Soln} can also be evaluated in $\mathcal O(1)$ operations resulting in a method to compute $\Phi(z)$ in $\mathcal O(1)$ operations.

\subsection{Fast evaluation of the derivative of generalized Hermite polynomials}

The function $U$ in~\eqref{e:Soln} satisfys the so-called zero-sum condition~\cite[Definition 3.5]{TrogdonSONNSD}.  This implies that differentiation commutes with the Cauchy integral operator $\CC_{\Gamma(n)}$. For each $j$, $U|_{\Gamma_{n,j}}'(z)$ is accurately computed with spectral differentiation.  Furthermore, $N'(z)$ can be computed accurately.  Define
\begin{align*}
T(z) = L_n^{-1}\tilde{\Phi}(z) = (\CC_{\Gamma(n)}[U](z)+ I)N(z) L_n^{-1} \begin{bmatrix} e^{n g(z)} & 0 \\ 0 & e^{-ng(z)} \end{bmatrix},
\end{align*}
and then
\begin{align*}
T'(z) &= \left(\CC_{\Gamma(n)}[U'](z)N(z) + (\CC_{\Gamma(n)}[U](z)+ I) N'(z) \phantom{\begin{bmatrix} 1 \\ 1 \end{bmatrix}}\right.\\
&\left.+  (\CC_{\Gamma(n)}[U](z)+ I) N'(z)\begin{bmatrix} ng'(z) & 0 \\ 0 & -ng'(z) \end{bmatrix} \right) L_n^{-1}\begin{bmatrix} e^{n g(z)} & 0 \\ 0 & e^{-ng(z)} \end{bmatrix}.
\end{align*}
Below, only the $(1,1)$ and $(2,1)$ entries of $T$ and $T'$ are needed.


\subsection{Newton's method for generalized Gauss--Hermite nodes}

In practice, we use Newton's method to find the zeros of
\begin{align*}
r_n(z) = T_{11}(z) e^{-V(z n^{1/(2m)})/2},
\end{align*}
which, of course, coincide with the zeros of $\pi_{n}(z)$.  Experiments show
that the magnitude of $r_n(z)$ is $\bigo(1)$ and hence it is more covenient for computation. For the $k$th node, one step of 
Newton's method takes the form:
\begin{align*}
\tilde x_k^{\rm new} = \tilde x_k^{\rm old} - \frac{T_{11}(\tilde x_k^{\rm old}) e^{-V(\tilde x_k^{\rm old} n^{1/(2m)})}}{[T_{11}'(\tilde x_k^{\rm old})- n^{1/(2m)} V'(\tilde x_k^{\rm old} n^{1/(2m)}) T_{11}(\tilde x_k^{\rm old})] e^{-V(\tilde x_k^{\rm old} n^{1/(2m)})}}.
\end{align*}
The exponential factors cancel out in this fraction but we leave them there as numerically we observe a 
small improvement in accuracy by doing so. If the update $|\tilde x_k^{\rm new} - \tilde x_k^{\rm old}|$ is sufficiently small, then 
we terminate Newton's method and the corresponding generalized Gauss--Hermite node is calculated via 
$x_k \approx \tilde x_k^{\rm new} n^{1/(2m)}$.

\subsection{Calculating generalized Gauss--Hermite weights}

Once the zeros $\{\tilde x_{k}\}_{k=1}^n$ of $\pi_n(x)$ are known the quadrature weights are found through the formula \cite{Hildebrand}
\begin{align}\label{weight-formula}
w_{k} = \frac{n^{1/(2m)}}{\gamma_{n-1} \pi_{n-1}(\tilde x_{k}) \pi_n'(\tilde x_{k})}    = - \frac{2 \pi n^{1/(2m)}}{T_{11}'(\tilde x_{k}) T_{21}(\tilde x_{k})}.
\end{align}
The monic polynomial $\pi_n$ has an exponentially small amplitude and it is convenient to choose a normalization so that the polynomial is typically $\bigo(1)$ on the support of $\psi_{V,n}$.  We have removed left multiplication by $L_n$ in the definition of $T(z)$ so that $T_{11}(z)$ is a more favorable multiple of $\pi_n$.

\subsection{Subsampling}

For each fix $n$, $V$, and $\epsilon > 0$, we can find a threshold parameter $\tau_{n,V,\epsilon}$ such that $|w_{k}| < \epsilon$ if $k < \tau_{n,V,\epsilon}$ or $k > n - \tau_{n,V,\epsilon}$ where $\epsilon$ is less than the smallest positive normalized floating-point number in double precision.  A node $x_{k}$ that satisfies $\tau_{n,V,\epsilon} \leq k \leq 1- \tau_{n,V,\epsilon}$ is said to be {\em non-trivial}.  If $V(x) = x^{2m} + \bigo(x^{2m -1})$, we demonstrate below that $\tau_{n,V,\epsilon} = \bigo(n^{1  - 1/(2m)})$ and give a method for choosing a constant $c$ so that $c n^{1  - 1/(2m)} \leq \tau_{n,V,\epsilon}$.

It follows from~\cite[(7.187), (7.84)]{Deift_Book}, assuming the differentiability of the asymptotic formula, that
\begin{align*}
\pi_{n-1}'(x) \pi_{n}(x) = (n C_n(x) + \bigo(1)) e^{-2 n \real g(x)}, \quad \gamma_{n-1} \sim \frac{1}{\pi 2^{3/2}} e^{n\ell},
\end{align*}
where $C_n(x) \leq C$ for all $n$.  It also follows that $e^{-2 n \real g(x)+n \ell} = e^{nV(x)}$, see~\cite[(7.49)]{Deift_Book}.  We assume that $1/C_n(x_{k})$ is bounded for all $k$ and $n$ and we have
\begin{align*}
w_k \leq \frac{D}{n} e^{-n V(x_{k})}.
\end{align*}
With the chosen scaling, $\psi_{V}(x) = \lim_{n \rightarrow \infty} \psi_{V,n}(x)$ gives the asymptotic density of the nodes $x_{k}$, which all lie in a finite interval that contains the origin.  Let $R_n>0$ be the largest value such that $\frac{D}{n} e^{-n V(\pm R_n)} \geq \epsilon$ and $R_n = \bigo(n^{-1/(2m)})$. Since $\psi_V$ is a continuous density we have
\begin{align*}
{\int_{-R_n}^{R_n} \psi_V(x) dx} = \bigo(n^{-1/(2m)}),
\end{align*}
which is an upper bound on the asymptotic fraction of non-trival nodes.  We find in practice that we can take $(n C_n(x) + \bigo(1))^{-1} \leq 1/n$ and if $V(x) = x^{2m}$ this gives
\begin{align*}
R_n \leq \frac{1}{n^{1/(2m)}} \left( \log{\epsilon^{-1}} + \log n^{-1} + \log (\pi 2^{3/2}) \right)^{1/(2m)}.
\end{align*}
Then
\begin{align*}
\tau_{n,V,\epsilon} \approx \left\lfloor n \int_{-\infty}^{-R_n} \psi_V(x) dx \right\rfloor,
\end{align*}
where the right-hand side is actually a lower bound for $\tau_{n,V,\epsilon}$ under the assumptions we have put forth.

\subsection{Examples}

\begin{figure}[tbp]
\centering
\begin{minipage}{.49\textwidth}
 \includegraphics[width=\textwidth]{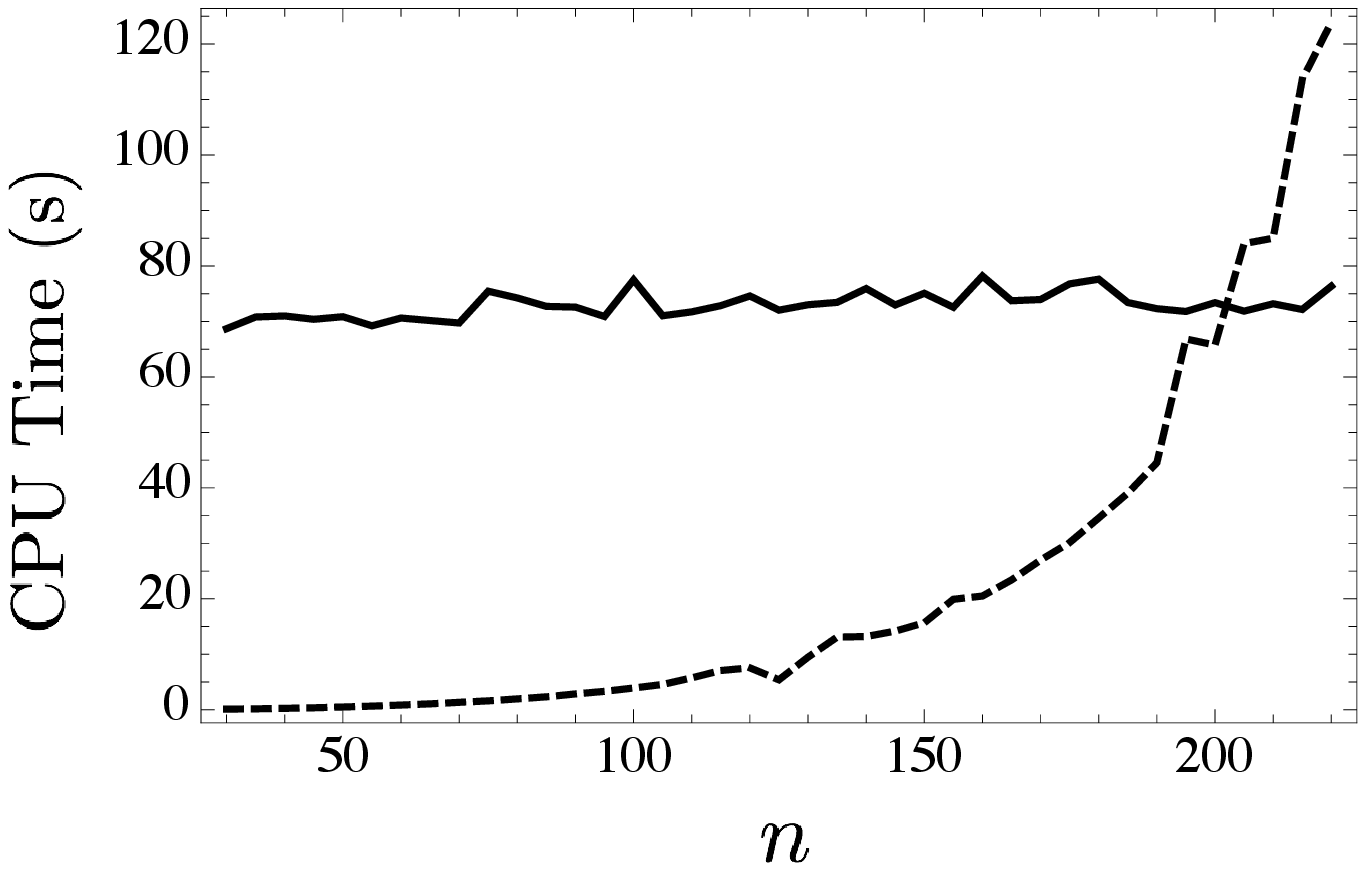}
\end{minipage}
\begin{minipage}{.49\textwidth}
\includegraphics[width=\textwidth]{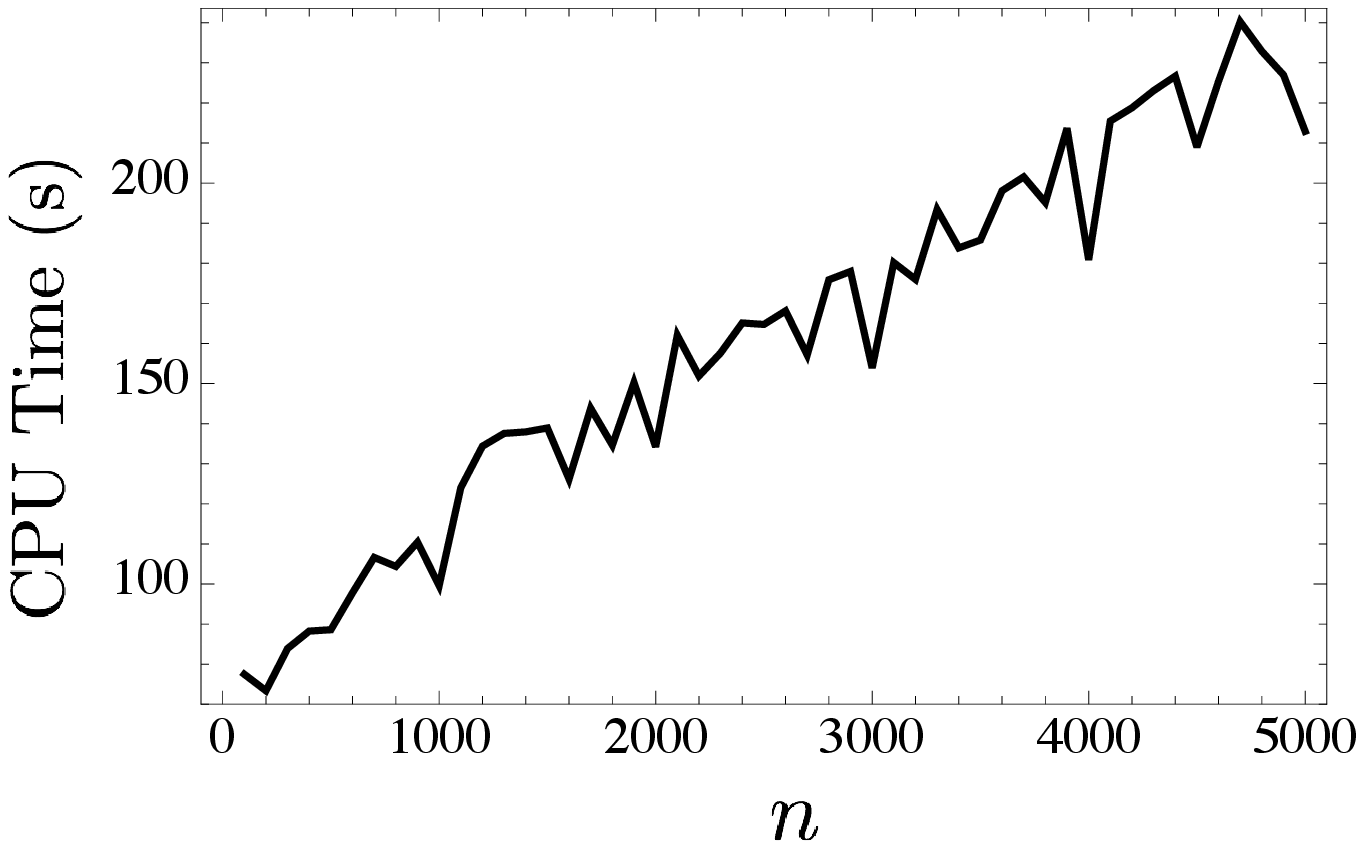}
\end{minipage}
\caption{Left: A comparison of the CPU time required to compute the recurrence coefficients with the Stieltjes procedure (dashed) with the CPU time required to compute the zeros of $\pi_n(x)$ with the approach advocated in this section (solid) when $V(x) = x^4$.  We see rapid growth with respect to $n$ for the Stieltjes procedure and linear growth for the Newton's method/RH approach.  Right: The CPU time required to compute the zeros of $\pi_n(x)$ with our approach for larger $n$.  There is clear linear growth. \label{RHCPUTime}}
\end{figure}

\begin{figure}[tbp]
\centering
\begin{minipage}{.49\textwidth}
 \includegraphics[width=\textwidth]{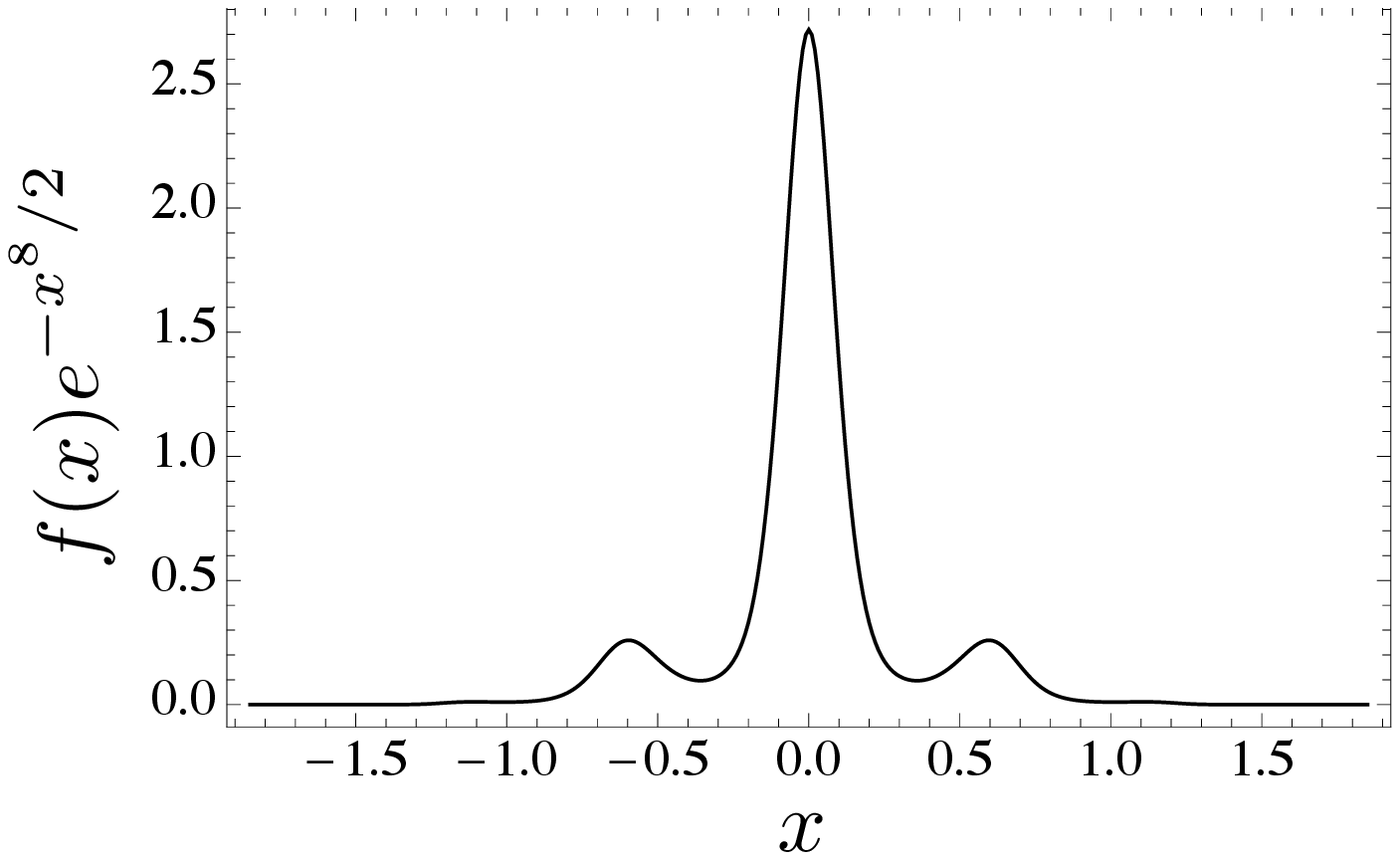}
\end{minipage}
\begin{minipage}{.49\textwidth}

\vspace{-.3cm}

\includegraphics[width=\textwidth]{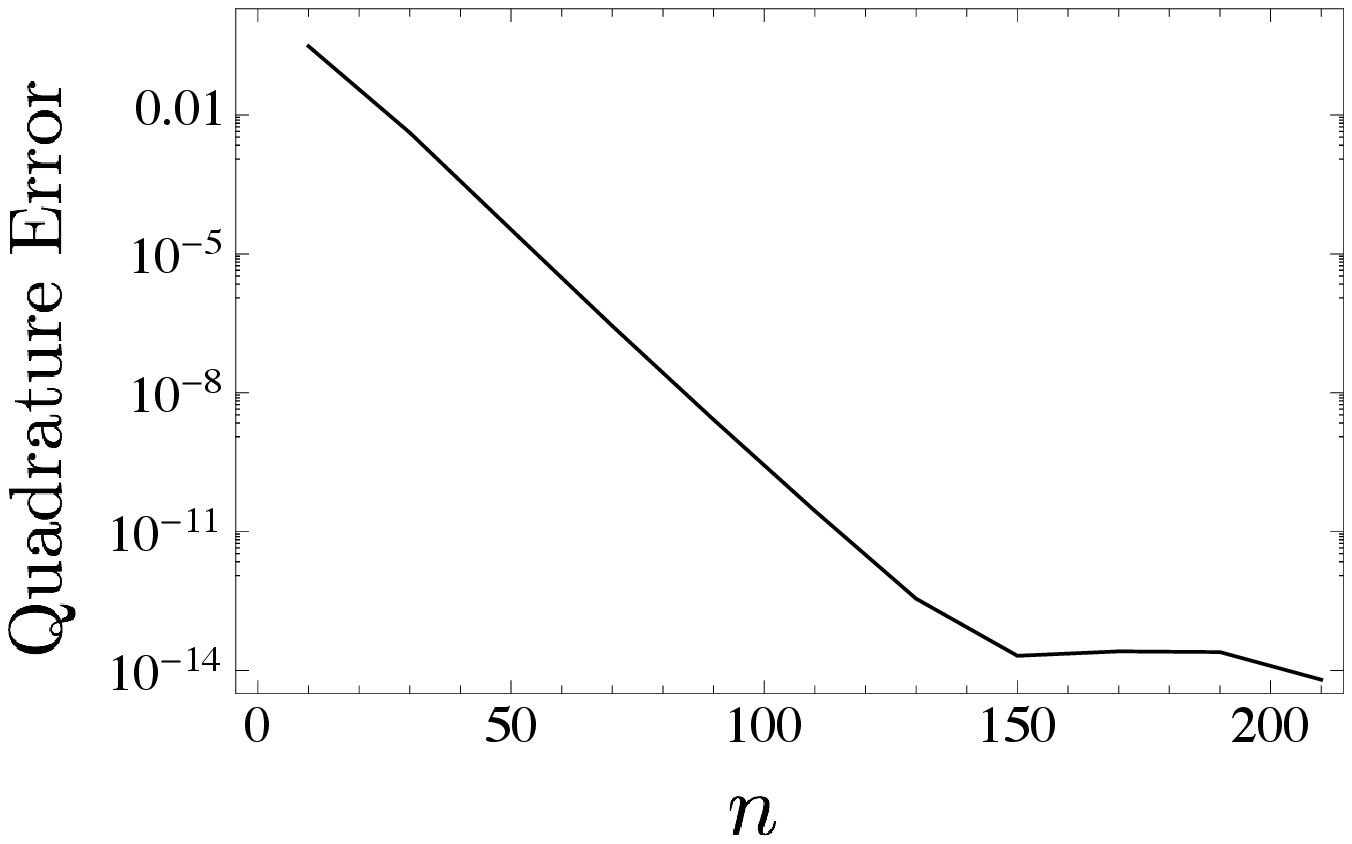}
 \end{minipage}
\caption{Left: A plot of $f(x)e^{-V(x)/2}$ for $f(x) = e^{\cos(10x)}/(1+25x^2)$ and $V(x) = x^8$. Right: The convergence 
of $\sum_{k=1}^n f(x_{k}) w_k$ to the integral of $f(x) e^{-V(x)}$.  The result is compared against Clenshaw--Curtis 
on $[-3,3]$ with 10,000 points. \label{EightInterpolationError1}}
\end{figure}

First, we compare the method described in this section with an implementation of the Stieltjes procedure because it appears to be the current method of choice to compute zeros when $V(x)$ is not quadratic \cite{Gautschi}.  The Stieltjes procedure is used to compute the coefficients in the $3$-term recurrence formula that the orthogonal polynomials satisfy.  These coefficients are used to construct an $n\times n$ tridiagonal matrix whose eigenvalues are the zeros of the $n$th-order polynomial.  For exponential weights, it appears that a straightforward implementation of this method has computational complexity that grows like $n^3$ \cite{TrogdonSOGaussQuad}.  In Figure~\ref{RHCPUTime} we demonstrate the efficiency of the method for large $n$ with $V(x) = x^4$.  We also note that our algorithm can be run for $V(x)=x^2$ to compare with the method above and errors on the order of $10^{-14}$ for the nodes are found.  We demonstrate quadrature errors in Figure~\ref{EightInterpolationError1} with $V(x) = x^8$.

\begin{figure}[tbp]
\centering
\includegraphics[width=.45\linewidth]{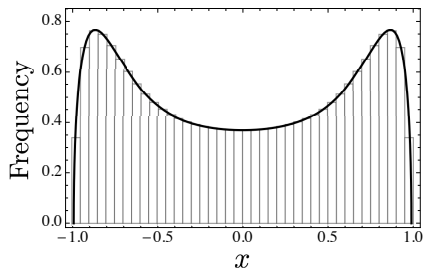}
\caption{A histogram for all of the approximated zeros $\{\tilde x_k\}$ of $\pi_{100000}(x)$ overlayed with the density $\psi_{V,n}(x)$ when $V(x) = x^8$.  The histogram shows that the distribution of the zeros approximates $\psi_{V,n}(x)$ in the sense of \eqref{dist-conv}. \label{EightHisto}}
\end{figure}

We can also examine the distribution of the zeros $\{\tilde x_k\}$ with a histogram.  In Figure~\ref{EightHisto} we show the equilibrium measure density $\psi_{V,n}(x)$ for $V(x) = x^8$ overlayed with a histrogram for $100,\! 000$  zeros.  Note that $\psi_{V,n}(x)$ is actually independent of $n$ because $V$ is a monomial.  It is clear from the figure that the distribution of the zeros approximates $\psi_{V,n}(x)$ in the sense of \eqref{dist-conv}.

\section{Application to Interpolation}\label{sec:interpolation}

\begin{figure}[tbp]
\centering
\includegraphics[width=.45\linewidth]{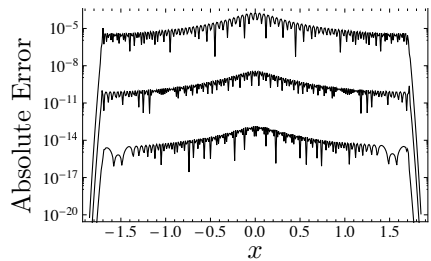}
\caption{The absolute error $|\tilde{\mathcal L}_n[f](x) - f(x)|e^{-V(x)/2}$ plotted as 
a function of $x$ for $n = 100, 200, 300$.  We see the expected spectral convergence of the approximation.  
This also validates the accuracy of the barycentric nodes and weights.\label{EightInterpolationError2}}
\end{figure}

Everything we have described here can be used for barycentric Lagrange interpolation.  We use the second form of the barycentric interpolation formula as discussed in \cite{BerrutTrefethen}:
\begin{align}\label{bc-form}
\mathcal L_n[f](x) = \sum_{k=1}^n \frac{f(x_{k}) \lambda_j}{x- x_{k}} / \sum_{k=1}^n \frac{\lambda_k}{x- x_k} , \quad \lambda_k = \frac{c}{\pi'_n(x_{k})},
\end{align}
for any convenient constant\footnote{Often, $c$ is chosen so that $\max_k |\lambda_k|=1$.} $c$. It is clear, in light of previous discussion, that we can evaluate $\{x_{k}\}_{k=1}^n$ and $\{\lambda_k\}_{k=1}^n$ in $\bigo(n)$ operations.   For reasons we discuss in Appendix~\ref{stability}, we evaluate
\begin{align}\label{weighted-bc}
\mathcal L_n[f](x)e^{-V(x)/2}. 
\end{align}
We demonstrate the convergence of an interpolant for $f(x) = e^{\cos(10x)}/(1+25x^2)$ with $V(x) = x^8$ in Figure~\ref{EightInterpolationError2}.

\section*{Conclusions}

We have shown the usefulness of a general algorithm for finding quadrature nodes and weights when the 
associated orthogonal polynomials and their derivatives can be pointwise evaluated in $\mathcal O(1)$ operations.  
The algorithm achieves optimal and state-of-the-art complexity of $\mathcal O(n)$ operations to compute $n$ 
quadrature nodes and weights.  In the case of Hermite polynomials the algorithm appears to achieve a smaller 
constant for this $\mathcal O(n)$ term when compared to other existing algorithms.  The method extends to 
quadrature with respect to general weights $e^{-V(x)}$ with the assistance of Riemann--Hilbert problems.  
We note that while $V(x)$ was chosen to be polynomial here, the extension can be made to entire weights 
(also known as Erd\H{o}s weights) such as $V(x) = \cosh x$~\cite{TrogdonSOGaussQuad}.

Presumably, extensions of these ideas can be made to the Laguerre weights $w(x) = x^{\alpha} e^{-V(x)}$ 
for $x \in (0,\infty)$.  Again, we expect the asymptotic expansion to be useful in the classical case 
of $V(x) = -x$ and the Riemann--Hilbert approach to aid in the general approach.  The Riemann--Hilbert 
approach may also apply to generalized Jacobi-type weights, \emph{i.e.}, $(1-x)^\alpha (1+x)^\beta e^{-V(x)}$ 
on $[-1,1]$, to extend the work of~\cite{Hale_13_01}.

The Riemann--Hilbert approach we employed used a numerical method to approximate the polynomials.  
A possible alternative would be to use asymptotic expansions of the polynomials derived from the Riemann--Hilbert 
approach.  An open question remains as to whether a high-order expansion can be computed effectively from 
the asymptotic Riemann--Hilbert theory.

\section*{Acknowledgements} 
We wish to thank Nick Trefethen for discussing this work with us and to 
Nick Hale for implementing the Glaser--Lui--Rokhlin algorithm in MATLAB.  We acknowledge the 
generous support of the National Science Foundation through grant NSF-DMS-130318 (TT). Any opinions, 
findings, and conclusions or recommendations expressed in this material are those of the authors and 
do not necessarily reflect the views of the funding sources. 

\appendix

\section{Weighted stability of Barycentric interpolation}\label{stability}

In this section we discuss the results of Higham~\cite{Higham} in the context of interpolation on $\mathbb R$.  Define the condition number for a function $f$ at a point $x$ by
\begin{align*}
\cond(f,x,n) = \lim_{\epsilon \rightarrow 0^+} \sup_{|\Delta f| \leq \epsilon |f|} \left\{ \left| \frac{\mathcal L_n[f](x) - \mathcal L_n[f+ \Delta f](x)}{\epsilon \mathcal L_n[f](x) } \right| \right\}.
\end{align*}
The main result of \cite{Higham} is the following theorem.

\begin{theorem}[\cite{Higham}]
If $u$ is machine unit roundoff and $\tilde {\mathcal L}_n[f](x)$ is the computed value, then
\begin{align*}
\left| \frac{\mathcal L_n[f](x) - \tilde{\mathcal L}_n[f](x)}{\mathcal L_n[f](x)} \right|  \leq (3n+4) u \cond(f,x,n) + (3n + 2) u \cond(1,x,n) + \bigo(u^2).
\end{align*}
\end{theorem}

This Theorem shows that if $\cond(1,x,n)$ and $\cond(f,x,n)$ do not grow too quickly with respect to $n$, then the barycentric formula \eqref{bc-form} is forward stable.  On the real line two questions remain:
\begin{enumerate}
\item How does $\cond(1,x,n)$ depend on $x$ and $n$?
\item For what class of functions $f$ is $\cond(f,x,n)$ relatively small?
\end{enumerate}
We begin by noting that $\mathcal L_n[\cdot]$ is a linear functional and $\mathcal L_n[1](x) = 1$ so that
\[
\cond(1,x,n) = \lim_{\epsilon \rightarrow 0^+} \sup_{|\Delta f| \leq \epsilon} \frac{|\mathcal L_n[\Delta f](x) |}{\epsilon} = \sup_{|\Delta f| \leq 1} |\mathcal L_n[\Delta f](x) | = \sum_{j=1}^n |\ell_{j,n}(x)|,
\]
where $\ell_{j,n}$ is the Lagrange polynomial of degree $n$ that takes the value of $1$ at $\tilde x_{j,n}$ and vanishes at all the other nodes.  The last equality can be seen by considering $\Delta f(\tilde x_{j,n}) = \sign \ell_{j,n}(x)$.  In the case of Chebyshev polynomials, this quantity grows like $\log n$, uniformly in $x \in [-1,1]$.  When the domain of interest is the whole real line we must 
modify the definition slightly as a quick numerical experiment shows that $\cond(1,x,n)$ grows exponentially.

One na\"\i ve bound gives
\begin{align*}
\sum_{j=1}^n |\ell_{j,n}(x)| e^{-V(x)/2}  \leq \sum_{j=1}^n |\ell_{j,n}(x)| e^{V(\tilde x_{j,n})/2-V(x)/2} =: \Lambda_n(x),
\end{align*}
where $\Lambda_n(x)$ is the so-called {\em weighted Lebesgue function} and $\Lambda_n = \sup_{x\in\mathbb{R}} \Lambda_n(x)$ is the {\em weighted Lebesgue constant}~\cite{LubinskyWeighted}.
For our choice of a exponential weighting it follows that $\Lambda_n = \bigo(n^{1/6})$ \cite{LubinskyWeighted}.

We are led to consider the following weighted relative error:
\begin{align*}
\begin{split}
\left| \frac{\mathcal L_n[f](x) - \tilde{\mathcal L}_n[f](x)}{\mathcal L_n[f](x)} \right|&e^{-V(x)/2}\\
\leq (3n&+4) u \cond(f,x,n)e^{-V(x)/2} + C n^{1/6} (3n + 2) u + \bigo(u^2).
\end{split}
\end{align*}
Therefore, we achieve forward stability when the errors are damped by the weight $e^{-V(x)/2}$.  We interpret this as meaning that the appropriate function space is \\$L^\infty(\mathbb R, e^{-V(x)/2}dx)$ as opposed to $L^\infty([-1,1])$.

Our final task is to identify a class of functions that is sufficiently regular so that $\cond(f,x,n)$ is well-behaved. 
By standard approximation theory we can relate the interpolation error to the best approximation error. 
\[
\|e^{-V(\cdot)/2} (f(\cdot) - \mathcal L_n[f](\cdot))\|_{L^\infty(\mathbb R)} \leq (1 + \Lambda_n ) \inf_{\deg P \leq n-1} \|e^{-V(\cdot)/2}(f(\cdot) - P(\cdot))\|_{L^\infty(\mathbb R)}.
\]

For a function that is absolutely continuous and of weighted bounded variation (see~\cite[Thm.~10.4]{LubinskyWeighted}), polynomial interpolation converges uniformly on compact sets. 

Assuming $n$ is sufficiently large so that $|\mathcal L_n[f](x)|^{-1} \leq 2 |f(x)|^{-1}$ we have
\begin{align*}
\cond(f,x,n) &\leq \frac{2}{|f(x)|} \sum_{j=1}^n |\ell_{j,n}(x)| \epsilon^{-1}|\Delta f(\tilde x_{j,n})| \\
&\leq  \frac{2}{|f(x)|} \sum_{j=1}^n |\ell_{j,n}(x)|e^{V(\tilde x_{j,n})/2} |\epsilon^{-1}\Delta f(\tilde x_{j,n})|e^{-V(\tilde x_{j,n})/2}\\
&\leq 2 \frac{\|e^{-V(\cdot)/2}f(\cdot)\|_{L^\infty(\mathbb R)}}{|f(x)|} \Lambda_n(x) e^{V(x)/2}
\end{align*}
so that $\cond(f,x,n)e^{-V(x)/2}$ behaves well with respect to both $x$ and $n$.  Not surprisingly, for functions that are well represented by polynomials the weighted barycentric formula \eqref{weighted-bc} is forward stable when errors are damped by the exponential function $e^{-V(x)/2}$.

\end{document}